\documentclass[11pt,reqno]{amsart}
\usepackage{amsmath,amsfonts}
\usepackage{amssymb}
\usepackage{amsthm}
\usepackage[all,cmtip,line]{xy}
\usepackage{array}
\usepackage{psfrag}
\usepackage{overpic}
\usepackage{bm}
\usepackage{setspace}

\newtheorem{theorem}{Theorem}
\newtheorem{lemma}{Lemma}[section]
\newtheorem{proposition}[lemma]{Proposition}
\newtheorem{remark}[lemma]{Remark}
\newtheorem{problem}{{\bf{Problem}}}

\theoremstyle{definition}

\newcommand \alp{\alpha}
\newcommand \eps{\varepsilon}
\newcommand \vphi{\varphi}

\newcommand \Gam{\Gamma}

\newcommand \gam{\gamma}
\newcommand \om{\omega}
\newcommand \tx{\text}

\newcommand \R{\mathbb{R}}

\newcommand \til{\tilde}

\newcommand \der{\partial}

\newcommand \mcl{\mathcal}

\newcommand \N{\mcl{N}}

\newcommand \ol{\overline}

\newcommand \Om{\Omega}
\newcommand \cd{\Gam_{D}}
\newcommand \gd{g_D}

\newcommand \NR{\mcl{N}_R}

\newcommand \NhRp{\mcl{N}_{R/2}^+}
\newcommand \NhRm{\mcl{N}_{R/2}^-}

\newcommand \Gamp{\Gam^+}
\newcommand \Gamm{\Gam^-}
\newcommand \Gamz{\Gam^0}

\newcommand \mt{m_0^+}
\newcommand \mb{m_0^-}

\newcommand \ult{u^+_{l}}
\newcommand \ulb{u^{-}_{l}}
\newcommand \urt{u^+_{r}}
\newcommand \urb{u^-_{r}}

\newcommand \pl{p_{l}}
\newcommand \pr{p_{r}}
\newcommand \St{S_0^+}
\newcommand \Sb{S_0^-}
\newcommand \Bt{B_0^+}
\newcommand \Bb{B_0^-}

\newcommand \rholt{\rho^+_{l}}
\newcommand \rholb{\rho^-_{l}}
\newcommand \rhort{\rho^+_{r}}
\newcommand \rhorb{\rho^-_{r}}
\newcommand \clt{c^+_l}
\newcommand \clb{c^-_l}

\newcommand \Ul{U_l}
\newcommand \Ur{U_r}
\newcommand \vphil{\vphi_l}
\newcommand \vphir{\vphi_r}

\newcommand \vphiz{\vphi_0}
\newcommand \Fz{\bm{F}_0}

\numberwithin{equation}{section}

\begin{document}
\title[Stability of contact discontinuity]
{Stability of contact discontinuity for steady Euler system in infinite duct}

\author{Myoungjean Bae}

\address{Department of Mathematics\\
         POSTECH\\
         San 31, Hyojadong, Namgu, Pohang, Gyungbuk, Republic of Korea}
\email{mjbae@postech.ac.kr or mybjean@gmail.com}


\date{}

\maketitle


\begin{abstract}
In this paper, we prove stability of contact discontinuities for full Euler system. We fix a flat duct $\N_0$ of infinite length in $\R^2$ with width $W_0$, and consider two uniform subsonic flow $\Ul^{\pm}=(u_l^{\pm}, 0,\pl,\rho_l^{\pm})$ with different horizontal velocity in $\N_0$ divided by a flat contact discontinuity $\Gam_{cd}$. And, we slightly perturb the boundary of $\N_0$ so that the width of the perturbed duct converges to $W_0+\omega$ for $|\omega|<\delta$ at $x=\infty$ for some $\delta>0$. Then, we prove that if the asymptotic state at left far field is given by $\Ul^{\pm}$, and if the perturbation of boundary of $\N_0$ and $\delta$ are sufficiently small, then there exists unique asymptotic state $\Ur^{\pm}$ with a flat contact discontinuity $\Gam_{cd}^*$ at right far field($x=\infty$) and unique weak solution $U$ of the Euler system so that $U$ consists of two subsonic flow with a contact discontinuity in between, and that $U$ converges to $\Ul^{\pm}$ and $\Ur^{\pm}$ at $x=-\infty$ and $x=\infty$ respectively. For that purpose, we establish piecewise $C^1$ estimate across a contact discontinuity of a weak solution to Euler system depending on the perturbation of $\der\N_0$ and $\delta$.

\end{abstract}

\bigskip

\section{Introduction}
Let $\rho, \bm u$ and $p$ be density, velocity and pressure of flow. Then
steady inviscid compressible flow is governed by the \emph{steady Euler system}
\begin{equation}
\label{1-a1}
\begin{split}
&div(\rho \bm u)=0\\
&div(\rho \bm u\otimes \bm u+p\mathbb I)=0\quad(\mathbb{I}:\tx{identity matrix})\\
&div(\rho \bm u B)=0
\end{split}
\end{equation}
with the \emph{Bernoulli's invariant} $B$ given by
\begin{equation}
\label{1-c1}
B=\frac 12|\bm u|^2+\frac{\gam p}{(\gam-1)\rho}
\end{equation}
for \emph{ideal polytropic gas} with an adiabatic exponent $\gam>1$. The quantity $c$, given by
\begin{equation}
\label{soundspeed}
c:=\sqrt{\frac{\gam p}{\rho}},
\end{equation}
is called the \emph{sound speed}. The flow type is classified by \emph{the Mach number} $M=\frac{|\bm u|}{c}$. If $M>1$ then the flow is called \emph{supersonic}, if $M<1$ then it is called \emph{subsonic}. If $M=1$ then the flow is called \emph{sonic}. If the flow is supersonic then the system \eqref{1-a1} is hyperbolic, if the flow is subsonic then the system \eqref{1-a1} becomes a hyperbolic-elliptic mixed system.

Due to the nonlinearity of the system \eqref{1-a1}, one expects that a solution of \eqref{1-a1} may contain discontinuities such as shocks or contact discontinuities even if a boundary condition is given by a smooth function. Such discontinuities can be described through a weak formulation of \eqref{1-a1}. A shock and a contact discontinuity are characterized by a normal velocity. While the normal velocity is nonzero on a shock, the normal velocity on a contact discontinuity completely vanishes. More details are given in Section \ref{Sec2}.
Because of the difference, one needs different schemes to study a shock and a contact discontinuity. For the case of a shock problem, one can identify a shock as a graph by using nonzero normal velocity on the shock. Owing to this advantage, the stability or instability of various shock phenomena have been investigated in many works(see \cite{Canic-Key-Lie}, \cite{Jch-Ch-F}, \cite{Ch-F3}, \cite{cs}, \cite{csf2}, \cite{CF}, \cite{kim}, \cite{TLiu2}, \cite{zh2} and references therein). For the case of a contact discontinuity, on the other hand, we need a different strategy due to the zero normal velocity on the contact discontinuity. For that reason, this subject has been studied in restricted regimes yet(\cite{cs}, \cite{csf}). In \cite{cs}, S. Chen proved the stability of steady Mach reflection configuration in a bounded region in $\R^2$ provided that an appropriate constant pressure is fixed on a cut-off boundary. In the Mach reflection configuration, two reflected shocks are separated by a contact discontinuity. In \cite{csf}, S. Chen and B. Fang proved the conditional stability of a reflection and a refraction of shocks occurred when an incident shock hits the interface, which is a contact discontinuity, of two different media. But still, the structural stability of a contact discontinuity under a general perturbation in unbounded domain is unknown. The main difficulty in study of a contact discontinuity is that the states on both sides of a contact discontinuity are unknown, so one needs to solve a free boundary problem with both sides of a free boundary to be determined. And, this is the main difference from a shock problem of Euler system.

The study of contact discontinuity is essential to understand a Mach reflection, which is one of important but difficult subject. When a vertical incident shock hits an inclined ramp, if the incident shock is relatively strong or the angle of the ramp is relatively small, then the incident shock is reflected at a point away from the boundary of the ramp, and two reflected shocks are formed at the reflection point with a contact discontinuity in between. This phenomenon is called \emph{Mach reflection,} named after Ernst Mach. Also, the Mach reflection for steady Euler system can be considered through shock polar analysis(see \cite{cs}). It is conjectured that the steady Mach reflection in $\R^2$ is structurally stable. In order to prove this conjecture, one needs to prove stability of a contact discontinuity along with two reflected shocks. In this paper, we prove structural stability of a contact discontinuity in an infinite duct where the flow on both sides of the contact discontinuity is subsonic. This is related to the case where two reflected shocks are transonic shocks in steady Mach reflection.

We fix a flat duct $\N_0$ of infinite length in $\R^2$ with width of $W_0$, and consider two uniform subsonic flow in $\N_0$ divided by a flat contact discontinuity $\Gam_{cd}$. Then we perturb the boundary of $\N_0$ with a small $C^{1,\alp}$ function so that the width of the perturbed duct converges to $W_0+\omega$ for $|\omega|<\delta$ at $x=\infty$ for some small constant $\delta>0$. Then we want to show that there exists two layers of subsonic flow divided by a contact discontinuity in the perturbed nozzle, and that the new contact discontinuity is a small perturbation of $\Gam_{cd}$.

It is a new feature that we allow for a perturbed contact discontinuity to converge to different asymptotic states at $x=\pm\infty$. Since the right asymptotic width $W_0+\omega$ of a perturbed nozzle is not necessarily same as the width $W_0$ at $x=-\infty$, we expect for the asymptotic pressure $p_{\infty}$ at $x=\infty$ to be different from the asymptotic pressure $p_{-\infty}$ at $x=-\infty$, and this yields two different asymptotic states for a perturbed contact discontinuity at $x=\pm\infty$. Details are given in Section \ref{Sec2}.

Another interesting aspect is application of the result from \cite{Li} to this work. Like a shock problem, a contact discontinuity problem is a free boundary problem. But unlike a shock problem, the states on both sides of a contact discontinuity $\Gam$ is unknown. Moreover since the normal velocity on $\Gam$ is zero, it is not clear how to locate a position of $\Gam$. So we use the Euler-Lagrange transformation to reformulate the contact discontinuity problem as a fixed boundary problem. Then the new difficulty is to find a weak solution of a first order nonlinear elliptic-hyperbolic mixed system so that the weak solution satisfies the Rankine-Hugoniot jump condition on a fixed boundary. If we can show that the weak solution is piecewise $C^1$ in two subregions divided by the flattened contact discontinuity $\til{\Gam}$, then a simple integration by parts shows that the weak solution indeed satisfies the Rankine-Hugoniot jump condition on $\til{\Gam}$. From this point of view, we employ the result of \cite{Li} to achieve piecewise $C^1$ regularity of weak solution for a nonlinear elliptic-hyperbolic mixed system. According to Theorem 1.1 of \cite{Li}, if a discontinuity boundary, which is a contact discontinuity in our case, is in $C^{1,\beta}$ for $0<\beta<1$, then the regularity of a corresponding weak solution to a uniformly elliptic equation is weaker than piecewise $C^{1,\beta}$ regularity. But, this may cause a difficulty in the iteration procedure which is the way to solve our main problem. Fortunately, the Euler-Lagrange transformation transforms a contact discontinuity to a flat boundary which is smooth. Hence, there is no deterioration  of regularity of weak solutions to elliptic equations in the iteration, and this is an advantage of the Euler-Lagrange transformation.

In Section \ref{Sec2}, we compute asymptotic states at far field in a perturbed duct, and use the Euler-Lagrange transformation to reformulate our problem as a fixed boundary problem. Then we state our main theorems.
In Section \ref{Sec3}, we establish piecewise $C^1$ regularity of weak solutions to Euler system, and use this estimate to prove the main theorems. The main difficulty would be uniform $L^{\infty}$ estimate of weak solutions to a class of uniformly elliptic equations in unbounded domain especially because we have two different asymptotic states at far field $x=-\infty$ and $x=\infty$.
\section{Problems and main theorems}
\label{Sec2}
\subsection{Asymptotic states at far field}
Let $\rho, \bm u$ and $p$ be the density, velocity and the pressure of flow respectively. Then
steady flow of compressible polytropic gas is governed by the  Euler system
\begin{align}
\label{1-6}
&div(\rho \bm u)=0\\
\label{1-7}
&div(\rho \bm u\otimes \bm u+p\mathbb I)=0\\
\label{1-8}
&div(\rho \bm u B)=0
\end{align}
with the Bernoulli's invariant
\begin{equation}
\label{B}
B=\frac 12|\bm u|^2+\frac{\gam p}{(\gam-1)\rho}
 \end{equation}
 for an adiabatic exponent $\gam>1$. We note that if $(\rho, \bm u, p)$ is in $C^1$, and satisfies \eqref{1-6}--\eqref{1-8}, then it also satisfies the transport equations
\begin{equation}
\label{2-b1}
\bm u\cdot \nabla(\frac{p}{\rho^{\gam}})=0\;\;\tx{and}\;\;\bm u\cdot \nabla B=0,
\end{equation}
and this means that the entropy($\frac{p}{\rho^{\gam}}$) and the Bernoulli's invariant $B$ are preserved along each streamline in $C^1$ flow.

We consider flow in $\R^2$. Let $u_1$ be the horizontal component of $\bm u$, and let $u_2$ be the vertical component of $\bm u$.
For an open and connected set $\Om\subset \R^2$, if $\bm U=(u_1,u_2,p,\rho)$ satisfies
\begin{equation}
\label{1}
\begin{split}
\int_{\Om}\rho \bm u\cdot D \xi
=\int_{\Om}(\rho u_k\bm u+p\hat{\bm e}_k)\cdot D\xi
=\int_{\Om}\rho\bm uB \cdot D \xi =0
\end{split}
\end{equation}
for any $\xi\in C_0^{\infty}(\Om)$ and $k=1,2$, then $\bm U$ is called a \emph{weak solution} to the Euler system \eqref{1-6}-\eqref{1-8} in $\Om$.
Suppose that $\Om$ is divided into two subsets $\Om^{\pm}$ by a non self-intersecting $C^1$ curve, and that
 $\bm U$ is $C^1$ in $\Om^{\pm}$ and $C^0$ in $\ol{\Om^{\pm}}$. Then one can easily check by integration by parts that $\bm U$ is a weak solution of the Euler system if and only if $\bm U$
satisfies \eqref{1-6}-\eqref{1-8} pointwisely in $\Om^{\pm}$ and the
\emph{Rankine-Hugoniot jump conditions}(abbreviated as \emph{R-H conditions} hereafter)
\begin{align}
\label{2-2}
&[\rho\bm u \cdot \bm n]_{\Gam}=0\\
\label{2-2-1}
&[\rho(\bm u\cdot \bm n)\bm u+p\bm n]_{\Gam}=\bm 0\\
\label{2-2-2}
&\rho \bm u\cdot\bm n[B]_{\Gam}=0
\end{align}
for a unit normal $\bm n$ on $\Gam$ where $[F]_{\Gam}$ is defined by
$[F(x)]_{\Gam}:
=F(x)|_{\overline{\Om^-}}-F(x)|_{\overline{\Om^+}}$ for $x\in \Gam.$
By \eqref{2-2}, the condition \eqref{2-2-1} can be rewritten as
\begin{equation}
\label{1-d}
\rho(\bm u\cdot \bm n)[\bm u\cdot\bm\tau]_{\Gam}=0,\quad [\rho(\bm u\cdot \bm n)^2+p]_{\Gam}=0
\end{equation}
where $\bm\tau$ denotes a unit tangential on $\Gam$.
For $\rho>0$ in $\Om$, the first condition in \eqref{1-d} implies either $\bm u\cdot \bm n=0$ on $\Gam$ or $[\bm u\cdot \bm\tau]_{\Gam}=0$. Suppose that $\bm U$ is discontinuous on $\Gam$. If $\bm u\cdot\bm n\neq 0$ and $[\bm u\cdot \bm\tau]_{\Gam}=0$ hold on $\Gam$, then  $\Gam$ is called a \emph{shock}. If $\bm u\cdot \bm n=0$ and $[\bm u\cdot \bm\tau]_{\Gam}\neq 0$, then $\Gam$ is called a \emph{contact discontinuity}. From this, we get the R-H conditions corresponding to a contact discontinuity as follows:
\begin{equation}
\label{1-a0}
\begin{split}
&\bm u\cdot \bm n=0\quad\text{on}\;\;\Gam,\quad [p]_{\Gam}=0.
\end{split}
\end{equation}

In $\R^2$, we consider a flat duct $\mcl{N}_0=\{(x,y)\in\R^2:x\in\R, -1<y<1\}$ of infinite length and two layers of uniform flow in $\N_0$ divided by the line $y=0$ with satisfying the following properties:
\begin{itemize}
\item[(i)] The velocity and density of top and bottom layers are given by positive constants $(\ult,0), \rholt$ and $(\ulb,0), \rholb$ respectively with $\ult\neq \ulb$;
\item[(ii)] The pressure of both top and bottom layers is given by a positive constant $\pl$;
\item[(iii)] The top and bottom layers are subsonic flow. In other words, there hold
\begin{equation}
\label{2-c3}
\frac{\ult}{\clt}<1\;\;\tx{and}\;\;\frac{\ulb}{\clb}<1\;\;\tx{for $\clt=\sqrt{\gam \pl/\rholt}$ and $\clb=\sqrt{\gam \pl/\rholb}$.}
 \end{equation}

\end{itemize}
From this, we define a piecewise constant vector $\Ul$ by
\begin{equation}
\label{def-bsl}
\Ul(x,y):=
\begin{cases}
(\ult,0,\pl,\rholt)=:\Ul^+&\text{for}\;\;y>0\\
(\ulb,0,\pl,\rholb)=:\Ul^-&\text{for}\;\;y<0.
\end{cases}
\end{equation}
Then, $\Ul$ is a weak solution of the Euler system with a contact discontinuity on the line $y=0$.

 Let $\eta$ be a smooth function satisfying
\begin{equation}
\label{1-g2}
\eta(x)=\begin{cases}
0&\tx{for}\;\;x\le -1\\
1&\tx{for}\;\;x\ge 1
\end{cases}, \quad 0\le\eta'(x)\le 10\;\;\tx{for all}\;x\in\R .
\end{equation}
Fix two constants $\omega_{\pm}$ with $|\omega_{\pm}|<1$, and let $h_{\pm}\in C^{1,\alp}(\R)$ be two functions satisfying
\begin{align}
\label{1-a}
&\|h_-+1+\om_-\eta\|_{H^1(\R)}+\|h_+-1-\om_+\eta\|_{H^1(\R)}\le \sigma,\\
\label{2-a4}
&\|h_-+1\|_{C^{1,\alp}(\R)}\le \sigma,\;\;\|h_+-1\|_{C^{1,\alp}(\R)}\le \sigma,\\
\label{2-a6}
&\lim_{R\to \infty}\|h_-+1+\omega_-\eta\|_{C^{1,\alp}(\R\setminus[-R,R])}=\lim_{R\to \infty}\|h_+-1-\omega_+\eta\|_{C^{1,\alp}(\R\setminus[-R,R])}=0
\end{align}
for small $\sigma>0$ to be determined later. To satisfy \eqref{2-a4} and \eqref{2-a6}, we assume that
\begin{equation}
\label{2-f2}
|\om_{\pm}|\le \sigma.
\end{equation}
For such functions $h_{\pm}$, let us set
\begin{equation*}
\til{\N}:=\{(x,y)\in \R^2: h_-(x)<y<h_+(x), x\in \R\}.
\end{equation*}
$\til{N}$ is a duct perturbed from $\mcl{N}_0$ by the functions $h_{\pm}$. Particularly, the width of $\til{N}$ at $x=\infty$ is changed to $2+(\omega_++\omega_-)$ from 2.
\begin{figure}[htp]
\centering
\begin{psfrags}
\psfrag{yt}[cc][][0.8][0]{$y=1\phantom{aaa}$}
\psfrag{yb}[cc][][0.8][0]{$y=-1\phantom{aaa}$}
\psfrag{ot}[cc][][0.8][0]{$\phantom{a}|\om_+|$}
\psfrag{ob}[cc][][0.8][0]{$\phantom{a}|\om_-|$}
\psfrag{N}[cc][][0.8][0]{$\til{\N}$}
\includegraphics[scale=.4]{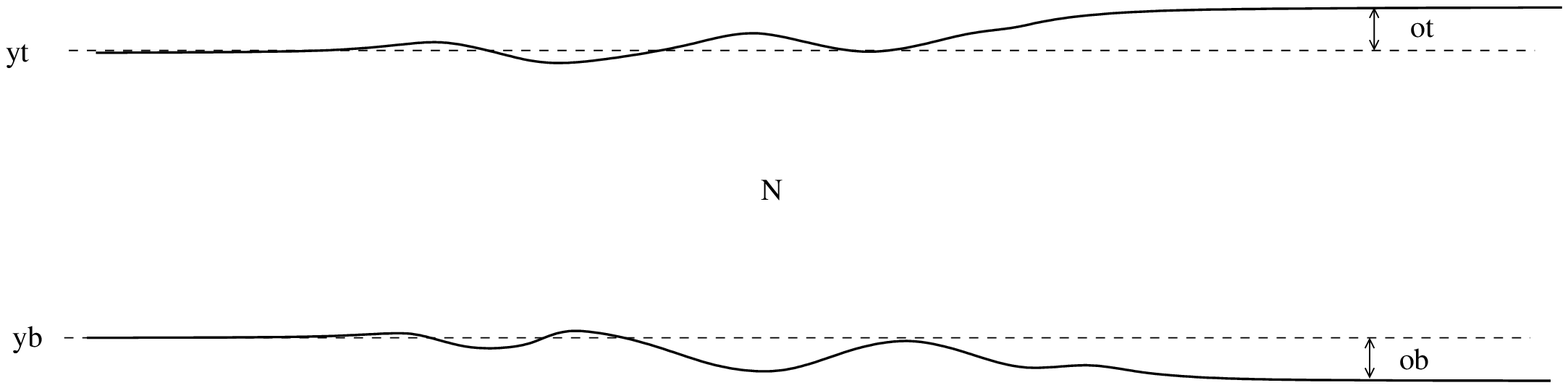}
\end{psfrags}
\end{figure}
The goal is to prove that there exists a weak solution of the Euler system in $\til{\N}$ with a contact discontinuity, and that the contact discontinuity is a small perturbation of $y=0$ for $\omega_{\pm}$ and $h_{\pm}$ satisfying \eqref{1-a}--\eqref{2-f2} where the parameter $\sigma$ is chosen sufficiently small. This is to achieve structural stability of contact discontinuity of steady Euler system even when the width of the duct is changed at far field. On the boundaries $y=h_{\pm}(x)$ of $\til{\N}$, we describe the slip boundary condition
\begin{equation}
\label{2-a7}
\bm u\cdot \bm n_{bd}=0
\end{equation}
where $\bm n_{bd}$ is the inward unit normal of $\til{N}$ on $y=h_{\pm}(x)$.

For a $C^1$ function $h_*(x)$ satisfying $h_-(x)<h_*(x)<h_+(x)$ for all $x$ in $\R$, suppose that
$$\
U(x,y)=(u, v, p, \rho)(x,y)
=\begin{cases}
(u+,v^+,p^+,\rho^+)(x,y)&\tx{for}\;\;y>h_*(x)\\
(u^-,v^-,p^-,\rho^-)(x,y)&\tx{for}\;\;y<h_*(x)
\end{cases}
$$
is a weak solution of the Euler system in $\til{N}$ with a contact discontinuity on $y=h_*(x)$ and $u>0$ in $\til{\N}$, and that $U$ is $C^1$ in $\til{\N}\setminus\{y=h_*(x)\}$, and satisfies the slip boundary condition \eqref{2-a7} on $\der\til{\N}$. Also, suppose that $U$ converges to $\Ul$ in $L^{\infty}$ at $x=-\infty$. Then the total mass flux on $x=x_0$ in each of top and bottom layers is preserved as well as the total mass flux on $x=x_0$ in $\til{N}$ is preserved for all $x_0$ in $\R$. In other words, $U$ satisfies
\begin{align}
\label{2-a8}
&\int_{h_-(x)}^{h_*(x)}(\rho^-u^-)(x,y)\;dy=\rholb\ulb=:m_0^-,\;
\int_{h_*(x)}^{h_+(x)}(\rho^+u^+)(x,y)\;dy=\rholt\ult=:m_0^+,\\
\label{2-a9}
&\int_{h_-(x)}^{h_+(x)}(\rho u)(x,y)\;dy=\rholb\ulb+\rholt\ult=:m_0.
\end{align}
\eqref{2-a8} and \eqref{2-a9} can be easily checked by using \eqref{1-6}, \eqref{1-a0} and \eqref{2-a7}. Moreover, by \eqref{2-b1} and positivity of $u$, $U$ also satisfies
\begin{equation}
\label{2-b2}
(S,B)(x,y)=\begin{cases}
\Bigl(\frac{\pl}{(\rholt)^{\gam}}, \frac 12 (\ult)^2+\frac{\gam \pl}{(\gam-1)\rholt}\Bigr)=:(\St,\Bt)&\tx{for}\;\;y>h_*(x)\\
\Bigl(\frac{\pl}{(\rholb)^{\gam}}, \frac 12 (\ulb)^2+\frac{\gam \pl}{(\gam-1)\rholb}\Bigr)=:(\Sb,\Bb)&\tx{for}\;\;y<h_*(x)
\end{cases}
\end{equation}
for $S=\frac{p}{\rho^{\gam}}$ and $B$ defined by \eqref{1-c1}.

By \eqref{2-a6}, the boundary of $\til{N}$ gets flat at far field. From this, we expect that the flow in $\til{N}$ at far field becomes two layers of uniform flow divided by a flat contact discontinuity. So we first compute the asymptotic state at far field in $\til{\N}$ consisting of two layers of uniform flow with a flat contact discontinuity in between. At $x=-\infty$, we fix $\Ul$ in \eqref{def-bsl} as the asymptotic state. It remains to compute the asymptotic state at $x=\infty$ corresponding to $\Ul$.
\begin{lemma}
\label{lemma-2-1}
Fix $\gam>1$, and let $\Ul$ be as in \eqref{def-bsl}. Set $\omega:=\om_++\om_-$ for $\om_{\pm}$ from \eqref{2-a6}. Then there is a constant $\om_0>0$ depending on $\Ul$ and $\gam$ such that for any $\omega\in[-\om_0,\om_0]$, there exist unique constants $\urt, \urb, \rhort, \rhorb$, $\pr$ and $\om_*\in(-1-\om_-, 1+\om_+)$ so that
\begin{equation}
\label{def-bsr}
\Ur(x,y)=\begin{cases}
(\urt,0,\pr,\rhort)=:\Ur^+&\tx{for}\;\;y>\om_*\\
(\urb,0,\pr,\rhorb)=:\Ur^-&\tx{for}\;\;y<\om_*
\end{cases}
\end{equation}
is a weak solution of \eqref{1-a1}, and satisfies \eqref{1-a0}\ and \eqref{2-a8}--\eqref{2-b2} with a contact discontinuity $y=\omega_*$. Furthermore, we have
$$
\frac{\der \pr}{\der \om}\ge \beta_0>0\;\;\tx{for all}\;\;\om\in[-\om_0,\om_0]
$$
where $\beta_0$ is a constant depending only on $\Ul,\gam$ and $\om_0$.
\begin{proof}
One can easily check that $\Ur$ in \eqref{def-bsr} is a weak solution of \eqref{1-a1}, and satisfies \eqref{1-a0}.
Set $\Ur(x,y)=(u_r, 0,\pr, \rho_r)(x,y)$.
By \eqref{2-b2}, we can write $(u_r,\rho_r)$ as follows:
\begin{equation*}
(u_r,\rho_r)(x,y)=\begin{cases}
\Bigl(\sqrt{2(\Bt-\frac{\gam (\St \pr^{\gam-1})^{1/\gam}}{(\gam-1)}}, (\frac{\pr}{\St})^{1/\gam}\Bigr)&\tx{for}\;\;y>\omega_*\\
\Bigl(\sqrt{2(\Bb-\frac{\gam (\Sb \pr^{\gam-1})^{1/\gam}}{(\gam-1)}}, (\frac{\pr}{\Sb})^{1/\gam}\Bigr)&\tx{for}\;\;y<\omega_*
\end{cases}.
\end{equation*}
From this, we define
\begin{equation*}
G(p, B,S):=(\frac pS)^{1/\gam}\sqrt{2\bigl(B-\frac{\gam}{\gam-1}(Sp^{\gam-1})^{1/\gam}\bigr)}.
\end{equation*}
Then, $\Ur$ satisfies \eqref{2-a8} and \eqref{2-a9} if and only if
\begin{align}
\label{2-c1}
&G(\pr,\Bt,\St)(1+\om_+-\om_*)=\mt,\\
\label{2-c2}
&G(\pr,\Bb,\Sb)(1+\om_-+\om_*)=\mb
\end{align}
for $m_0^{\pm}$ given by \eqref{2-a8}.
Solving \eqref{2-c1} for $\om_*$, and plugging it into \eqref{2-c2}, we get
\begin{align}
\label{2-e7}
&\om_*=1+\om_+-\frac{\mt}{G(\pr, \Bt,\St)},\\
\label{2-e8}
&H(\pr, \om):=G(\pr, \Bb,\Sb)(2+\om-\frac{\mt}{G(\pr,\Bt,\St)})-\mb=0\;\;\tx{for}\;\;\om=\om_++\om_-.
\end{align}
From \eqref{def-bsl} and \eqref{2-a8}, we have $H(\pl, 0)=0$. Also, a direct computation using \eqref{B} and \eqref{2-c3} yields
\begin{equation*}
\der_pG(\pl, B_0^{\pm},S_0^{\pm})=\frac{\pl^{\frac{1}{\gam}-1}((u_l^{\pm})^2-(c_l^{\pm})^2)}{\gam(S_0^{\pm})^{1/\gam}\sqrt{2(B_0^{\pm}-\frac{\gam (S_0^{\pm})^{1/\gam}\pl^{1/\gam}}{\gam-1})}}<0,
\end{equation*}
and this implies $\der_pH(\pl,0)<0$. By the implicit function theorem, we can choose a constant $\omega_0>0$ small depending on $\Ul$ so that for any $\omega\in[-\omega_0,\omega_0]$, there exists unique $\pr(\om)$ satisfying
$H(\pr(\om),\omega)=0$, and such $\pr(\om)$ is $C^1$ with respect to $\om$. Moreover, we may adjust $\omega_0$ to satisfy
\begin{equation*}
\frac{\der\pr}{\der\om}=-\frac{\der_w H}{\der_pH}>0\;\;\tx{for all}\;\;\om\in[-\om_0,\om_0].
\end{equation*}
Once $\pr(\om)$ is obtained, then $\om_*$ is given by \eqref{2-e7}.
\end{proof}
\end{lemma}
\begin{remark}
From Lemma \ref{lemma-2-1}, we have that if $\om_{\pm}$ satisfy \eqref{2-f2} with $\sigma\le \frac{\om_0}{2}$ for $\om_0$ in Lemma \ref{lemma-2-1}, then
\begin{equation}
\label{2-f3}
|\Ur^+-\Ul^+|+|\Ur^--\Ul^-|+|\om_*|\le C\sigma
\end{equation}
for a constant $C$ depending only on $\Ul$ and $\gam$.
\end{remark}

\subsection{Problems and main theorems}
The following is the main problem of this paper.
\begin{problem}
\label{problem1}
Fix $\om_{\pm}$ satisfying \eqref{2-f2} and $h_{\pm}$ satisfying \eqref{1-a}--\eqref{2-a6} for $\sigma\in(0,\om_0]$. Then, find a function $\gd$ and a vector valued function $U=(u, v, p, \rho)$ in $\til{\N}$ to satisfy the following properties:
\begin{itemize}
\item[(i)] $U$ is a weak solution to \eqref{1-6}-\eqref{1-8} in $\til{\N}$ in sense of \eqref{1};
\item[(ii)] $U$ is in $C^1$ in $\til{\N}\setminus \cd$ for
$
\cd:=\{(x,\gd(x)):x\in \R\};
$
\item[(iii)] $U$ satisfies the R-H conditions \eqref{1-a0} on $\cd$;
\item[(iv)] $\frac{|\bm u|}{c}<1$ in $\til{\N}$ for $|\bm u|=\sqrt{u^2+v^2}$ and $c$ defined by \eqref{soundspeed};
\item[(v)] $U$ satisfies the slip boundary condition \eqref{2-a7} on $\der\til{\N}$;

\item[(vi)]
Setting $\til{\N}^+:=\til{\N}\cap\{y>\gd(x)\}\quad\tx{and}\quad\til{\N}^-:=\til{\N}\cap\{y<\gd(x)\},$
\begin{equation}
 \label{2-b6}
 \begin{split}
&\lim_{R\to \infty}\|U-\Ul^+\|_{L^{\infty}(\til{\N}^+\setminus \{x\ge -R\})}=\lim_{R\to \infty}\|U-\Ul^-\|_{L^{\infty}(\til{\N}^-\setminus \{x\ge -R\})}=0,\\
 &\lim_{R\to \infty}\|U-\Ur^+\|_{L^{\infty}(\til{\N}^+\setminus \{x\le R\})}=\lim_{R\to \infty}\|U-\Ur^-\|_{L^{\infty}(\til{\N}^-\setminus \{x\le R\})}=0
 \end{split}
 \end{equation}
 where $\Ur$ is as in Lemma \ref{lemma-2-1};
 \item[(vii)] $u>0$ and $\rho>0$ hold in $\N$.
\end{itemize}
\end{problem}
In order to find a solution to Problem \ref{problem1}, we use weighted H\"{o}lder norms. For a connected open set $\Om\subset \R^2$, let $\Upsilon$ be a closed portion of the boundary of $\Om$. For $\bm x=(x,y), \bm {\til x}=(\til x,\til y)\in \Om$, set
\begin{equation*}
\begin{split}
&d_{\bm x}:=dist(\bm x, \Upsilon),\quad d_{\bm x, \bm {\til x}}:=\min(d_{\bm x}, d_{\bm {\til x}}).\\
\end{split}
\end{equation*}
For $k\in \R, \alp\in(0,1)$ and $m\in \mathbb{Z}_+$, we define
\begin{align*}
&\|u\|_{m,0,\Om}:=\sum_{0\le|\beta|\le m}\sup_{\bm x\in \Om}|D^{\beta}u(\bm x)|,\quad
[u]_{m,\alp,\Om}:=\sum_{|\beta|=m}\sup_{\bm x, \bm y\in\Om,\bm x\neq \bm y}\frac{|D^{\beta}u(\bm x)-D^{\beta}u(\bm y)|}{|\bm x-\bm y|^{\alp}}\\
&\|u\|_{m,0,\Om}^{(k,\Upsilon)}:=\sum_{0\le|\beta|\le m}\sup_{\bm x\in \Om}(d_{\bm x})^{\max(|\beta|+k,0)}|D^{\beta}u(\bm x)|\\
&[u]_{m,\alp,\Om}^{(k,\Upsilon)}:=\sum_{|\beta|=m}\sup_{\bm x, \bm {\til x}\in\Om,\bm x\neq \bm {\til x}}(d_{\bm x,\bm {\til x}})^{\max(m+\alp+k,0)}\frac{|D^{\beta}u(\bm x)-D^{\beta}u(\bm {\til x})|}{|\bm x-\bm {\til x}|^{\alp}}\\
&\|u\|_{m,\alp,\Om}:=\|u\|_{m,0,\Om}+[u]_{m,\alp,\Om},\quad \|u\|_{m,\alp,\Om}^{(k,\Upsilon)}:=\|u\|_{m,0,\Om}^{(k,\Upsilon)}+[u]_{m,\alp,\Om}^{(k,\Upsilon)}\\
\end{align*}
where we write $D^{\beta}=\der_{x}^{\beta_1}\der_{y}^{\beta_2}$ for a multi-index $\beta=(\beta_1,\beta_2)$ with $\beta_j\in\mathbb{Z}_+$ and $|\beta|=\beta_1+\beta_2$. $C^{m,\alp}_{(k,\Upsilon)}(\Om)$ denotes completion of the set of all smooth functions whose $\|\cdot\|_{m,\alp,\Om}^{(k,\Upsilon)}$-norms are finite in the norm $\|\cdot\|_{m,\alp,\Om}^{(k,\Upsilon)}$. Similarly, $C^{m,\alp}(\ol{\Om})$ denotes completion of the set of all smooth functions whose $\|\cdot\|_{m,\alp,\Om}$-norms are finite in the corresponding norm. Hence, $C^{m,\alp}_{(k,\Upsilon)}(\Om)$ and $C^{m,\alp}(\ol{\Om})$ are Banach spaces.
For a vector valued function $\bm Q=(q_j)_{j=1}^{N}$, let us set
\begin{equation*}
\|\bm Q\|_{m,\alp,\Om}:=\sum_{j=1}^N\|q_j\|_{m,\alp,\Om},\quad
\|\bm Q\|_{m,\alp,\Om}^{(k, \Upsilon)}:=\sum_{j=1}^N\|q_j\|_{m,\alp,\Om}^{(k, \Upsilon)}.
\end{equation*}
The following theorem indicates that if $\sigma$ is sufficiently small, then Problem \ref{problem1} has a solution.
\begin{theorem}
\label{theorem-1}Fix $\alp\in(0, 1)$, and fix $\om_{\pm}$ and $h_{\pm}$ with satisfying \eqref{2-f2} and \eqref{1-a}--\eqref{2-a6}. Then there exist constants $C_0>0$ and $\sigma_0\in(0,\om_0]$ depending on $\Ul, \gam$ and $\alp$ so that for any $\sigma\in(0,\sigma_0]$, there exists unique solution $U=(u,v,p,\rho)$ of Problem \ref{problem1} with a contact discontinuity $y=\gd(x)$ satisfying the following estimates:
\begin{itemize}
\item[(i)]
$\|\gd\|_{1,\alp, \R}\le C_0\sigma$;
\item[(ii)]
$\underset{R\to\infty}{\lim} \|\gd-\omega_*\eta\|_{1,\alp,\R\setminus[-R, R]}=0
$ for $\eta$ defined by \eqref{1-g2};
\item[(iii)] For $\til{\N}^{\pm}$ as in Problem \ref{problem1}(vi),
$U$ satisfies
\begin{equation*}
\|U-\Ul^+\|_{1,\alp,\til{\N}^+}^{(-\alp, \Gam_{h_+}\cup \cd)}+
\|U-\Ul^-\|_{1,\alp,\til{\N}^-}^{(-\alp, \Gam_{h_-}\cup \cd)}\le C_0\sigma;
\end{equation*}
\item[(iv)] $U$ converges to the asymptotic states $\Ul$ at $x=-\infty$ and $\Ur$ at $x=\infty$ in the following sense:
\begin{equation*}
\begin{split}
&\lim_{R\to \infty}\|U-U_l^+\|_{1,\alp, \til{\N}^+\setminus\{x\ge -R\}}^{(-\alp, \Gam_{h_+}\cup \cd)}=
\lim_{R\to \infty}\|U-U_l^-\|_{1,\alp, \til{\N}^-\setminus\{x\ge -R\}}^{(-\alp, \Gam_{h_-}\cup \cd)}=0,\\
&\lim_{R\to \infty}\|U-U_r^+\|_{1,\alp, \til{\N}^+\setminus\{x\le R\}}^{(-\alp, \Gam_{h_+}\cup \cd)}=
\lim_{R\to \infty}\|U-U_r^-\|_{1,\alp, \til{\N}^-\setminus\{x\le R\}}^{(-\alp, \Gam_{h_-}\cup \cd)}=0;
\end{split}
\end{equation*}
\item[(v)] $\rho\ge \frac 1{10}\min(\rholt, \rholb)>0$ and $u\ge \frac 1{10}\min(\ult, \ulb)>0$ hold in $\til{\N}$.
\end{itemize}
\end{theorem}

\begin{theorem}
\label{theorem-1-2}
Fix $\om_{\pm}$ and $h_{\pm}$ with satisfying \eqref{2-f2} and \eqref{1-a}--\eqref{2-a6} with $\sigma\in(0,\sigma_0]$ for $\sigma_0$ in Theorem \ref{theorem-1}. Then, there exists a constant $C_1>0$ depending only on $\Ul$ and $\gam$ so that for any $(U,\gd)$ satisfying all the properties in Theorem \ref{theorem-1}, there hold
\begin{align}
\label{1-g4}
&\|U-U_0\|_{L^2(\til{\N})}\le C_1\sigma,\\
\label{1-g5}
&\|\gd-\omega_*\eta\|_{L^2(\R)}\le C_1\sigma
\end{align}
for $U_0(x,y)=\begin{cases}
(1-\eta(x))(\ult, 0,\pl,\rholt)+\eta(x)(\urt, 0,\pr,\rhort)&\tx{for}\;\;y>\gd(x)\\
(1-\eta(x))(\ulb, 0,\pl,\rholb)+\eta(x)(\urb, 0,\pr,\rhorb)&\tx{for}\;\;y<\gd(x)
\end{cases}$.
\end{theorem}

From \eqref{1-a0}, the tangential of $\cd$ is parallel to the velocity of flow on both sides of $\cd$. So if we use the \emph{Euler-Lagrange coordinate transformation}, then $\cd$ becomes a fixed flat boundary in the new coordinates while $\cd$ is a free boundary to be found simultaneously with a weak solution $U$ of the Euler system to solve Problem \ref{problem1}. Moreover, due to the conservation of total mass flux in $x$-direction in $\til{\N}^{\pm}$, $\Gam_{h_{\pm}}$ become flat in the new coordinates.

Let $(\gd, U)$ be a solution of Problem \ref{problem1}, and for $(x,y)\in \til{\N}$, define a transformation $T$ by
\begin{equation}
\label{1-1}
T=(T_1, T_2): (x,y)\mapsto (x,\int_{h_-(x)}^y \rho u(x,t)\;dt-m_0^-).
\end{equation}
By \eqref{1-6}, \eqref{2-a7}, Problem \ref{problem1}(ii),(v) and \eqref{1-a0}, we have
\begin{equation}
\label{1-c}
\begin{split}
\frac{d}{dx}T_2(x,\gd(x))=\frac{d}{dx}T_2(x,h_{\pm}(x))=0\;\;\text{for any}\;\; x\in \R.
\end{split}
\end{equation}
By \eqref{1-c} and \eqref{2-b6}, we obtain
\begin{equation}
\label{1-a4}
\begin{split}
&T(\til{\N})=\R\times(-m^-_0, m^+_0)=:\N,\\
&T(\cd)=\{(X,Y)\in \R^2: Y=0\}=:\Gamz,\\
&T(\Gam_{h_-})=\{(X,Y)\in \R^2: Y=-m_0^-\}=:\Gamm,\\
&T(\Gam_{h_+})=\{(X,Y)\in \R^2: Y=m_0^+\}=:\Gamp.
\end{split}
\end{equation}
For convenience, we use the notations of
\begin{equation*}
\begin{split}
&\N^+=\N\cap\{(X,Y)\in \R^2: Y>0\},\quad \N^-=\N\cap \{(X,Y)\in \R^2: Y<0\}.
\end{split}
\end{equation*}
Since $\det DT=\rho u$, $T$ is invertible if $U$ satisfies Theorem \ref{theorem-1}(v). Then $V=(u,v,p,\rho)(X,Y)$ given by $V=U\circ T^{-1}$ is well defined, and it becomes a weak solution of the following system in $\N$:
\begin{align}
\label{1-5}
&\der_X(\frac{1}{\rho u})-\der_Y(\frac vu)=0\\
\label{1-10}
&\der_X(u+\frac{p}{\rho u})-\der_Y(\frac{pv}{u})=0\\
\label{1-11}
&\der_X v+\der_Y p=0\\
\label{1-12}
&\frac 12(u^2+v^2)+\frac{\gam p}{(\gam-1)\rho}=\begin{cases}
\Bt&\tx{for}\;\;Y>0\\
\Bb&\tx{for}\;\;Y<0
\end{cases}(=:B_0).
\end{align}
If a weak solution $V$ of \eqref{1-5}--\eqref{1-12} is in $C^1(\N^{\pm})\cap C^0(\ol{\N^{\pm}})$, then the integration by parts combined with the fact that $\nu=(0,1)$ is a unit normal  on $\Gamz$ yields the R-H condition
\begin{equation}
\label{1-a8}
[\frac vu]_{\Gamz}=[p]_{\Gamz}=0.
\end{equation}
Conversely, we can show that if $V\in C^1(\N^{\pm})\cap C^0(\ol{\N^{\pm}})$ is a weak solution of \eqref{1-5}--\eqref{1-12} then there exist inverse Euler-Lagrange transformation $\mathfrak{T}$ so that $U=V\circ \mathfrak{T}^{-1}$ satisfies all the properties of Problem \ref{problem1}.
\begin{lemma}
\label{lemma-1}
Fix two $C^{1,\alp}$ functions $h_{\pm}$.
There exists a constant $\delta_1>0$ depending only on $\Ul$ and $\gam$ so that if $V=(u,v,p,\rho)$ satisfies the following properties:
\begin{itemize}
\item[(i)] $\|V-\Ul^+\|_{L^{\infty}(\N^+)}+\|V-\Ul^-\|_{L^{\infty}(\N^-)}\le \delta_1$;
\item[(ii)] $V$ is a weak solution of \eqref{1-5}-\eqref{1-12}, and $V$ is in $C^1(\N^{\pm})\cap C^0(\ol{\N^{\pm}})$;
\item[(iii)] $V$ satisfies the slip boundary condition
\begin{equation}
\label{2-b8}
\frac{v}{u}(X,-m_0^-)=h_-'(X)\;\;\text{on}\;\;\Gamm,\quad \frac vu(X, m_0^+)=h_+'(X)\;\;\text{on}\;\;\Gamp
\end{equation}
\end{itemize}
then,
\begin{itemize}
\item[(a)] the transformation
\begin{equation}
\label{1-15}
\mathfrak{T}:(X,Y)\mapsto (X, \int_{-m_0^-}^Y\frac{1}{\rho u}(X,t)dt+h_-(X))=:(x,y)
\end{equation}
is well defined, and satisfies $\mathfrak{T}=T^{-1}$ for $T$ defined by \eqref{1-1};
\item[(b)] $U=V\circ \mathfrak{T}^{-1}$ satisfies (i)-(v) of Problem \ref{problem1} with $\cd$ given by
\begin{equation*}
\cd=T^{-1}(\{Y=0\})=\{(x,y): y= \int_{-m_0^-}^0\frac{1}{\rho u}(x,t)dt+h_-(x)\}.
\end{equation*}
\end{itemize}
\begin{proof}
By using \eqref{1-5}, \eqref{2-b8}, \eqref{1-15}, and choosing $\delta_1$ small, one can directly check that $D\mathfrak{T}=(DT)^{-1}$. Also, we have $\mathfrak T(\{Y=-m_0^-\})=\{y=h_-(x)\}$. This proves (a).

 Using \eqref{1-5}--\eqref{1-12} and \eqref{2-b8}, one can directly show that $\bm U$ satisfies properties (i),(ii),(iv),(v) of Problem \ref{problem1} if we reduce $\delta_1>0$ depending only on $\Ul$ and $\gam$. By the R-H condition \eqref{1-a8}, $U$ satisfies
$
\rho({\bm u}\cdot \nu)[{\bm u}]_{\cd}=0
$
for $\bm u=(u,v)$
where $\nu$ indicates a unit normal on $\cd$,
and this implies either ${\bm u}\cdot\nu=0$ or $[{\bm u}]_{\cd}=0$ on $\cd$. From the definition of $\Ul$ in \eqref{def-bsl}, we have $ \ult\neq \ulb$. we further reduce $\delta_1$ to satisfy $\delta_1<\frac{|\ult-\ulb|}{10}$. Then we get $|[\bm u]_{\cd}|>0$. So we have $\bm u\cdot\nu=0$ on $\cd$.
\end{proof}
\end{lemma}
By Lemma \ref{lemma-1}, Problem \ref{problem1} can be replaced by the following problem:
\begin{problem}
\label{problem3}
Fix $\om_{\pm}$ satisfying \eqref{2-f2} and $h_{\pm}$ satisfying \eqref{1-a}--\eqref{2-a6} for $\sigma\in(0,\om_0]$ where $\om_0$ is as in Lemma \ref{lemma-2-1}. Then, find a vector valued function $V=(u, v, p, \rho)$ in $\N$ satisfying the following properties:
\begin{itemize}
\item[(i)] $V\in C^1(\N^{\pm})\cap C^0(\ol{\N^{\pm}})$ is a weak solution of \eqref{1-5}-\eqref{1-11}, and satisfies \eqref{1-12} in $\N$;
\item[(ii)] $V$ is subsonic in $\N$ in sense that $(u^2+v^2)< c^2$ for $c$ defined by \eqref{soundspeed};
\item[(iii)] $V$ satisfies the slip boundary condition \eqref{2-b8};
\item[(iv)] $V$ converges to the asymptotic states $\Ul$ at $x=-\infty$ and $\Ur$ at $x=\infty$ in the following sense:
\begin{equation*}
\begin{split}
&\lim_{R\to \infty}\|V-\Ul^+\|_{L^{\infty}(\N^+\setminus\{X\ge -R\})}=\lim_{R\to \infty}\|V- \Ul^-\|_{L^{\infty}(\N^-\setminus\{X\ge -R\})}=0,\\
&\lim_{R\to \infty}\|V-\Ur^+\|_{L^{\infty}(\N^+\setminus\{X\le R\})}=\lim_{R\to \infty}\|V- \Ur^-\|_{L^{\infty}(\N^-\setminus\{X\le R\})}=0;
\end{split}
\end{equation*}
\item[(v)] $u>0$ and $\rho>0$ hold in $\N$.
\end{itemize}
\end{problem}
Hereafter, we write $(X,Y)$ as $(x,y)$. If $V=(u, v, p, \rho)$ is in $C^1$, then using \eqref{1-5}, \eqref{1-11} and \eqref{1-12}, we can rewrite \eqref{1-10} as
\begin{equation*}
\der_x(\frac{p}{\rho^{\gam}})=0.
\end{equation*}
So if $V$ is a solution to Problem \ref{problem3}, then we have
\begin{equation}
\label{2-b9}
\frac{p}{\rho^{\gam}}=\begin{cases}
\St&\;\;\text{for}\;\;y>0\\
\Sb&\;\;\text{for}\;\;y<0
\end{cases}
(=:S_0).
\end{equation}
From \eqref{1-5}, we expect that there is a function $\vphi$ satisfying
\begin{equation}
\label{4-5}
\vphi_x=\frac vu,\quad  \vphi_y=\frac{1}{\rho u}
\end{equation}
 so that we can rewrite \eqref{1-11} and \eqref{1-12} in terms of $\vphi_x, \vphi_y$, $S_0$ and $B_0$ as in \cite{Jch-Ch-F}. If so, by \eqref{2-b8}, $\vphi$ should satisfy
\begin{equation}
\label{k}
\vphi=h_-+k_-\;\;\text{on}\;\;\Gamm,\quad \vphi=h_++k_+\;\;\text{on}\;\;\Gamp
\end{equation}
for some constants $k_{\pm}$. First, we define $\vphi$ corresponding to the asymptotic states $\Ul$ and $\Ur$ as follows:
\begin{equation}
\label{2-f1}
\vphil(x,y)=\begin{cases}
\frac{y}{m_0^+}(=:\vphil^+)&\tx{for}\;\;y>0\\
\frac{y}{m_0^-}(=:\vphil^-)&\tx{for}\;\;y<0
\end{cases},\;\; \vphir(x,y)=\begin{cases}
\frac{y}{\rhort\urt}+\omega_*(=:\vphir^+)&\tx{for}\;\;y>0\\
\frac{y}{\rhorb\urb}+\omega_*(=:\vphir^-)&\tx{for}\;\;y<0
\end{cases}.
\end{equation}
Then, we have
\begin{equation*}
\vphil(x,y)=\lim_{x\to-\infty}h_{\pm}(x)\;\;\tx{and}\;\;\vphir(x,y)=\lim_{x\to \infty}h_{\pm}(x)\;\;\tx{on}\;\;\Gam^{\pm}.
\end{equation*}
From this, we choose $k_{\pm}=0$ in the boundary condition \eqref{k}.

From \eqref{2-b9}, we have $p=S_0\rho^{\gam}$, and by plugging this into $(\rho^2\cdot \eqref{1-12})$ and using \eqref{4-5}, we get
\begin{equation}
\label{2-e1}
\mathfrak G(\rho, D\vphi)=0
\end{equation}
for $\mathfrak G$ defined by
\begin{equation}
\label{2-e6}
\mathfrak G(\rho, \bm q):=B_0\rho^2-\frac{\gam}{\gam-1}S_0\rho^{\gam+1}-\frac{1+q_1^2}{2q_2^2}\;\;\text{for}\;\;\bm q=(q_1, q_2)\in \R^2.
\end{equation}
A direct computation shows $
\mathfrak G(\rho_l^{\pm}, D\vphil^{\pm})=0$ and $\der_{\rho}\mathfrak G(\rho_l^{\pm}, \bm D\vphil^{\pm})=-\rho_l^{\pm}((c_l^{\pm})^2-(u_l^{\pm})^2)<0\;\;\text{in}\;\;\N^{\pm},
$
so we can choose a constant $\delta_2>0$ depending only on $\Ul$ and $\gam$ so that
for any $(x,y,\bm q)\in (\N^+\times B_{\delta_2}(D\vphil^+))\cup (\N^-\times B_{\delta_2}(D\vphil^-))$, there exists unique $\rho=\rho(x,y,\bm q)$ satisfying
\begin{equation}
\label{3-e7}
\mathfrak{G}(\rho(x,y,\bm q), \bm q)=0,
\end{equation}
and such $\rho(x,y,\bm q)$ is continuously differentiable with respect to $\bm q$. We write as $\rho(x,y,\bm q)$ rather than $\rho(\bm q)$ because $B_0$ in \eqref{1-12} and $S_0$ in \eqref{2-b9} are piecewise constant functions.

Let $\vphi\in C^0(\ol{\N})\cap C^1(\ol{\N^{\pm}})$ be a function satisfying
$
\|\vphi-\vphil^+\|_{C^1(\ol{\N^+})}+\|\vphi-\vphil^-\|_{C^1(\ol{\N^-})}\le \delta_2,
$
then there exists unique $\rho(x,y, D\vphi)$ satisfying the equation
\begin{equation}
\label{2-e3}
\mathfrak{G}(\rho(x,y,D\vphi), D\vphi)=0\;\;\text{in}\;\;\N^{\pm}.
\end{equation}
For such $\rho(x,y,D\vphi)$, we use \eqref{2-b9} and \eqref{4-5} to express $u, v$ and $p$ as
\begin{equation}
\label{1-b3}
 u=\frac{1}{\rho\vphi_y},\quad v=\frac{\vphi_x}{\rho\vphi_y},\quad p=S_0\rho^{\gam}
\end{equation}
so that the equation
\eqref{1-11} can be rewritten as
\begin{equation}
\label{4-6}
div\bigl(\bm A(x,y,D\vphi)\bigr)=0\;\;\text{in}\;\;\N
\end{equation}
where $\bm A(x,y,\bm q)=(A_1, A_2)(x,y,\bm q)$ is given by
\begin{equation}
\label{4-7}
A_1(x,y,\bm q)=\frac{q_1}{\rho(x,y,\bm q)q_2},\quad A_2(x,y,\bm q)=S_0(x,y) \rho^{\gam}(x,y,\bm q)
\end{equation}
with
\begin{equation}
\label{2-a3}
\|A_j\|_{C^k(\ol{\N^+\times B_{{3\delta_2}/{4}}(D\vphil^+)})}+\|A_j\|_{C^k(\ol{\N^-\times B_{{3\delta_2}/{4}}(D\vphil^-)})}\le C_k
\end{equation}
for each $k\in \mathbb N$ and $j=1,2$, where the constant $C_k$ depends only on $\Ul, \gam$ and $k$.

Next, we consider the R-H condition for $\vphi$ on $\Gamz$. Rewriting \eqref{1-a8} in terms of $\vphi$, we get $[\vphi_x]_{\Gamz}=[A_2(x,y, D\vphi)]_{\Gamz}=0$. In particular, we may rewrite $[\vphi_x]_{\Gamz}=0$ as $[\vphi]_{\Gamz}=k_*$ for a constant $k_*$. Furthermore, we choose $k_*=0$ by continuity of $\vphil$ and $\vphir$ across $\Gamz$. Because, we will seek a solution $\vphi$ of \eqref{4-6} so that $\vphi$ converges to $\vphil$ at $x=-\infty$ and to $\vphir$ at $x=\infty$. So we get
\begin{equation}
\label{1-16}
[\vphi]_{\Gamz}=[A_2(x,y, D\vphi)]_{\Gamz}=0.
\end{equation}
Now we consider the boundary value problem \eqref{4-6}, \eqref{k} with $k_{\pm}=0$ and \eqref{1-16}.

\begin{theorem}
\label{theorem-2}
Fix $\alp\in(0, 1)$, and fix $\om_{\pm}$ satisfying \eqref{2-f2} and $h_{\pm}$ satisfying \eqref{1-a}--\eqref{2-a6}. Then, there are constants $C_2$ and $\sigma_1\in(0,\om_0]$ depending on $\Ul, \gam$ and $\alp$ so that wherever $\sigma\in(0,\sigma_1]$, the boundary value problem \eqref{4-6} and \eqref{k} with $k_{\pm}=0$ has unique weak solution $\vphi\in H^1_{loc}(\N)$ satisfying the following properties:
\begin{itemize}
\item[(i)] $\vphi$ is in $C^0(\ol{\N})\cap C^1(\ol{\N^{\pm}})\cap C^2(\N^{\pm})$, hence $\vphi$ satisfies \eqref{4-6} in $\N^{\pm}$, and \eqref{1-16} on $\Gamz$ pointwisely;
\item[(ii)] The equation \eqref{4-6} is uniformly elliptic in $\N$;
\item[(iii)] $\vphi$ satisfies the estimate
\begin{equation}
\label{2-e4}
\|\vphi-\vphil^+\|_{2,\alp,\N^+}^{(-1-\alp, \Gamp\cup\Gamz)}+\|\vphi-\vphil^-\|_{2,\alp,\N^-}^{(-1-\alp, \Gamm\cup\Gamz)}\le C_2\sigma;
\end{equation}
\item[(iv)] $\vphi$ converges to $\vphil$ at $x=-\infty$ and to $\vphir$ at $x=\infty$ in the following sense:
\begin{equation}
\label{2-e5}
\begin{split}
&\lim_{R\to \infty}\|\vphi-\vphil^+\|_{2,\alp,\N^+\setminus\{x\ge -R\}}^{(-1-\alp, \Gamp\cup\Gamz)}=
\lim_{R\to \infty}\|\vphi-\vphil^-\|_{2,\alp,\N^-\setminus\{x\ge -R\}}^{(-1-\alp, \Gamm\cup\Gamz)}=0,\\
&\lim_{R\to \infty}\|\vphi-\vphir^+\|_{2,\alp,\N^+\setminus\{x\le R\}}^{(-1-\alp, \Gamp\cup\Gamz)}=
\lim_{R\to \infty}\|\vphi-\vphir^-\|_{2,\alp,\N^-\setminus\{x\le R\}}^{(-1-\alp, \Gamm\cup\Gamz)}=0;
\end{split}
\end{equation}
\item[(v)] $\der_y\vphi\ge \frac{1}{m_*}$ holds for some $m_*>0$ in $\N$ where $m_*>0$ depends only on $\Ul$ and $\gam$.
\end{itemize}
\end{theorem}

We first prove Theorem \ref{theorem-2}, then prove Theorem \ref{theorem-1} and Theorem \ref{theorem-1-2}.

\section{Proof of  Theorem \ref{theorem-2}}
\label{Sec3}

In order to prove Theorem \ref{theorem-2}, we need to prove unique existence of a weak solution $\vphi\in H^1_{loc}(\N)$ to
the boundary value problem
\begin{align}
\label{3-a1}
&div\bigl(\bm A(x,y, D\vphi)\bigr)=0\quad\text{in}\;\;\N,\\
\label{3-a2}
&\vphi=h_-\quad\text{on}\;\;\Gamm,\quad \vphi=h_+\quad\text{on}\;\;\Gamp
\end{align}
for $\bm A$ defined by \eqref{4-7}. Furthermore, the weak solution $\vphi$ is required to satisfy the additional R-H condition \eqref{1-16} on $\Gamz$. If the equation \eqref{3-a1} is strictly elliptic and $\vphi\in C^1(\ol{\N^{\pm}})\cap C^2(\N^{\pm})$, then the weak Harnack inequality implies that $\vphi$ is continuous across $\Gamz$, and $[A_2(x,y,D\vphi)]_{\Gamz}=0$ easily follows from the integration by parts because the vector $(0,1)$ is unit normal of $\Gamz$. Therefore the key point to prove Theorem \ref{theorem-2} is to show that the boundary value problem \eqref{3-a1}, \eqref{3-a2} has a piecewise $C^1$ weak solution. For that purpose, we employ results from \cite{Li}.

We prove Theorem \ref{theorem-2} in two steps. First, we formulate a linearized boundary value problem where coefficients of an elliptic equation in the boundary problem are piecewise $C^{\alp}$, and apply the result of \cite{Li} to weak solutions of the boundary value problem. Main difficulty in the first step would be uniform $L^{\infty}$ estimate of weak solutions in unbounded domain $\N$ because of two different asymptotic states $\vphil$ and $\vphir$ at $x=\pm\infty$. Then, we use a fixed point theorem to prove Theorem \ref{theorem-2}.

\subsection{Linearized boundary value problem}
\label{subsec-3-1}
Define a smooth connection from $\vphil$ to $\vphir$ as follows:
For $\eta$ defined by \eqref{1-g2}, we set
\begin{equation}
\label{3-b7}
\vphiz(x,y)=(1-\eta(x))\vphil(x,y)+\eta(x)\vphir(x,y)
\end{equation}
for $\vphil$ and $\vphir$ defined by \eqref{2-f1}. By \eqref{2-f3}, if we choose $\sigma$ small, then $\bm A(x,y,D\vphiz)$ is well defined by \eqref{2-a3}.  Since $div(\bm A(x,y,D\vphil))=0$, \eqref{3-a1} is equivalent to
$
div(\bm A(x,y,D\vphi)-\bm A(x,y,D\vphiz))=div \bm \Fz\;\;\tx{in}\;\;\N
$ for
\begin{equation}
\label{F}
\bm \Fz:=\bm A(x,y,D\vphil)-\bm A(x,y,D\vphiz).
\end{equation}
To ensure that $\Fz$ is well defined, we let
\begin{equation}
\label{1-g7}
\sigma\le \min\{\om_0, \frac{3\delta_2}{4C}\}=:\sigma^{\sharp}
\end{equation}
for $C$ in \eqref{2-f3}.
From \eqref{2-f1}, \eqref{1-b3}, \eqref{4-7},\eqref{3-b7} and the definition of $\Fz$ in \eqref{F}, we easily get
\begin{equation}
\label{F1}
\Fz(x,y)=\begin{cases}
(0,0)&\tx{for}\;\;x\le -1\\
(0, \pl-\pr)&\tx{for}\;\;x\ge 1.
\end{cases}
\end{equation}
Moreover, \eqref{2-f3} implies that $\Fz\in C^{\infty}(\ol{\N^+})\cap C^{\infty}(\ol{\N^-})$ satisfies the estimate
\begin{equation}
\label{F2}
\|\Fz\|_{C^k(\ol{\N^+})}+\|\Fz\|_{C^k(\ol{\N^-})}\le C_k\sigma
\end{equation}
for a constant $C_k$ depending only on $\Ul, \gam$ and $k$ for each $k\in \mathbb{Z}_+$.
For a fixed function $\phi$, set
\begin{equation}
\label{2-a}
a_{ij}^{(\phi)}(x,y):=\int_0^1D_{q_j}A_{i}(x,y,D\vphi_0+tD(\phi-\vphi_0))\;dt
\end{equation}
for $A_i(x,y,\bm q)$ defined by \eqref{4-7}.
$\vphi$ solves the boundary value problem of \eqref{3-a1} and \eqref{3-a2} if and only if $\psi:=\vphi-\vphi_0$ solves
\begin{equation}
\label{3-a4r}
\begin{split}
&\sum_{i,j=1}^2\der_i(a^{(\vphi)}_{ij}\der_j\psi)=div \Fz\quad\text{in}\;\;\N,\\
&\psi=h_-+1+(1-\eta)\omega_-=:g_-\quad\text{on}\;\;\Gamm,\quad \psi=h_+-1-(1-\eta)\om_+=:g_+\quad\text{on}\;\;\Gamp
\end{split}
\end{equation}
where $\der_1$ and $\der_2$ denote $\frac{\der}{\der x}$ and $\frac{\der}{\der y}$ respectively.

Fix $\alp\in(0,1)$, and define an iteration set $\mcl{K}_M$ by
\begin{equation}
\label{3-a3}
\mcl{K}_M:=\{\vphi\in C^{1,\alp}(\ol{\N^{\pm}})\cap C^0(\ol{\N}):\|\vphi-\vphil^+\|^{(-1-\alp, \Gamp\cup\Gamz)}_{2,\alp,\N^+}+\|\vphi-\vphil^-\|^{(-1-\alp, \Gamm\cup\Gamz)}_{2,\alp,\N^-}\le M\sigma\}
\end{equation}
for constants $M>1$ and $\sigma\in(0,\frac{\sigma^{\sharp}}{2}]$ to be determined later with $M\sigma\le \frac{3\delta_2}{4}$ for $\delta_2$ in \eqref{2-a3}. Then, $\mcl{K}_M$ is a convex and compact subset of the  Banach space $C^{1,\alp/2}_{(-1)}(\ol{\N^+})\cap C^{1,\alp/2}_{(-1)}(\ol{\N^-})$(see  Section 5.1 of \cite{Ch-F3} for details).
For a fixed $\phi\in \mcl{K}_M$, consider the following linear boundary value problem:
\begin{align}
\label{3-a4}
&\sum_{i,j=1}^2\der_i(a^{(\phi)}_{ij}\der_j\psi)=div \Fz\quad\text{in}\;\;\N,\\
\label{3-a5}
&\psi=h_-+1+(1-\eta)\omega_-=:g_-\quad\text{on}\;\;\Gamm,\quad \psi=h_+-1-(1-\eta)\om_+=:g_+\quad\text{on}\;\;\Gamp.
\end{align}
The following lemma is essential to prove Theorem \ref{theorem-2}.
\begin{proposition}
\label{proposition-1}Let $\delta_2$ be as in \eqref{2-a3}.
Fix $\alp\in(0, 1)$. Then there exist constants $C>0$ and $\eps_0\in(0,\frac{3\delta_2}{4}]$ depending only on $\Ul, \gam$ and $\alp$ so that if $M\sigma\le \eps_0$, then for any $\phi\in \mcl{K}_M$, the boundary value problem of \eqref{3-a4} and \eqref{3-a5} has unique weak solution $\psi\in H^1_{loc}(\N)\cap C^{0}(\ol{\N})$ with satisfying
\begin{align}
\label{3-a6}
&\|\psi\|_{2,\alp,\N^+}^{(-1-\alp, \Gamp\cup\Gamz)}+\|\psi\|_{2,\alp,\N^-}^{(-1-\alp, \Gamm\cup\Gamz)}\\
&\phantom{aaaa}
\le C(\|\Fz\|_{1,\alp,\N^+}+\|\Fz\|_{1,\alp,\N^-}+\|g\|_{1,\alp,\N}+\|Dg\|_{L^2(\N)}),\notag\\
\label{3-a7}
\tx{and}\;\;&\lim_{R\to \infty}\|\psi\|_{2,\alp,\N^+\setminus\{|x|\le R\}}^{(-1-\alp, \Gamp\cup\Gamz)}=
\lim_{R\to \infty}\|\psi\|_{2,\alp,\N^-\setminus\{|x|\le R\}}^{(-1-\alp, \Gamm\cup\Gamz)}=0
\end{align}
where $g$ is defined by \eqref{3-f4}.
\end{proposition}
In order to prove Proposition \ref{proposition-1}, we set a linear boundary value problem in a bounded domain as follows:
Let $\chi_0$ be a function satisfying
\begin{align*}
&\chi_0(y)=\begin{cases}
1&\text{for}\;\;y\le -\frac{m_0^-}{2},\\
0&\text{for}\;\;y\ge \frac{m_0^+}{2}
\end{cases},\quad 0\le \chi_0(y)\le 1\;\;\text{for}\;\;y\in[-\frac{m_0^-}{2}, \frac{m_0^+}{2}],\\
&\chi'_0(y)\le 0\;\;\text{for all}\;\;y\in\R\;\;\tx{and}\;\;\|\chi_0\|_{C^2(\R)}\le  C(m_0^{\pm})
\end{align*}
for a constant $C(m_0^{\pm})$ depending only on $m_0^{\pm}$ for $m_0^{\pm}$ in \eqref{2-a8}, and define a function $g$ by
 \begin{equation}
 \label{3-f4}
 g(x,y):=g_-(x)\chi_0(y)+g_+(x)(1-\chi_0(y))
 \end{equation}
 for $g_{\pm}$ given by \eqref{3-a4r}. Then $g$ satisfies the estimate
\begin{equation*}
\|g\|_{1,\alp,\N}\le C(\|g_-\|_{1,\alp,\R}+\|g_+\|_{1,\alp,\R})
\end{equation*}
for a constant $C$ depending on $\Ul$ and $\alp$. For a fixed constant $R\ge 10$, set
$
\widehat{\N}_R:=\N\cap\{(x,y)\in \R^2:|x|<R+2\},
$
and let $\N_R$ be a convex and connected domain satisfying
$
\widehat{N}_R\subset \N_R\subset \widehat{N}_{R+1}
$ where $\der\N_R$ is a simple closed smooth curve.
We consider the following boundary value problem in a bounded domain $\N_R$:
\begin{align}
\label{2-a1}
&\sum_{i,j=1}^2\der_i(a^{(\phi)}_{ij}\der_j\psi)=div \Fz\quad\text{in}\;\;\N_R,\\
\label{2-a2}
&\psi=g\quad\text{on}\;\;\der\N_R.
\end{align}
\eqref{2-a1} and \eqref{2-a2} has unique weak solution $\psi_R$ in $H^1(\N_R)$.
We claim that $\psi=\underset{R\to\infty}{\lim}\psi_R$ is unique weak solution of \eqref{3-a4} and \eqref{3-a5} with satisfying \eqref{3-a6} and \eqref{3-a7}.

\begin{lemma}
\label{lemma-2}Let $\delta_2$ be as in \eqref{2-a3}.
There exist positive constants $C>0$, $\eps_1\in(0,\frac{3\delta_2}{4}]$ and $\lambda$ depending on $\Ul$ and $\gam$ with $C$ depending on $\alp$ in addition so that if $M\sigma\le \eps_1$ then, for any $\phi\in \mcl{K}_M$, the coefficient matrix $\bigl(a_{ij}^{(\phi)}\bigr)_{i,j=1}^2$ defined by \eqref{2-a} satisfies the following properties:
\begin{itemize}
\item[(i)]
$
\lambda|\bm{\xi}|^2\le \sum_{i,j=1}^2a^{(\phi)}_{ij}(x,y)\xi_i\xi_j\le \frac{1}{\lambda}|\bm{\xi}|^2\quad\text{in}\;\;\N
$
for any $\bm{\xi}=(\xi_1, \xi_2)\in \R^2$;
\item[(ii)]$\bigl(a_{ij}^{(\phi)}\bigr)_{i,j=1}^2$ is symmetric, that is,
$
a_{12}^{(\phi)}(x,y)=a_{21}^{(\phi)}(x,y)\quad\text{in}\;\;\N;
$
\item[(iii)]
$
\|a_{ij}^{(\phi)}-a_{ij}^{(\vphil)}\|_{1,\alp,\N^+}^{(-\alp,\Gamp\cup\Gamz)}+\|a_{ij}^{(\phi)}-a_{ij}^{(\vphil)}\|_{1,\alp,\N^-}^{(-\alp,\Gamm\cup\Gamz)}
\le CM\sigma
$
for $i,j=1,2$.

\end{itemize}
\begin{proof}
Set $c^2=\gam S_0\rho^{\gam-1}$. Then a direct computation using \eqref{2-e6} and \eqref{3-e7} yields
$
\der_{q_1}\rho=\frac{- vu\rho}{c^2-|\bm u|^2}$ and $\der_{q_2}\rho=\frac{\rho^2u(u^2+v^2)}{c^2-|\bm u|^2}$ for
$|\bm u|^2=u^2+v^2$, and this implies
\begin{equation}
\label{1-b8}
\begin{split}
&\der_{q_1}A_1=u\frac{c^2-u^2}{c^2-|\bm u|^2}, \;\;
\der_{q_2}A_2=\frac{\rho^2c^2u|\bm u|^2}{c^2-|\bm u|^2},\quad \der_{q_2}A_1=\der_{q_1}A_2=-\frac{\rho c^2uv }{c^2-|\bm u|^2}.
\end{split}
\end{equation}
This proves (ii).
Using \eqref{2-a} and \eqref{1-b8}, we can easily show that
$
a^{(\vphil)}_{11}a^{(\vphil)}_{22}-\bigl(a^{(\vphil)}_{12}\bigr)^2
=\frac{\rho_l^2u_l^4c_l^2}{c_l^2-u_l^2}>0
$
for
$
(\rho_l, u_l, p_l)=\begin{cases}
(\rholt, \ult, \pl)&\;\;\text{for}\;\;y>0\\
(\rholb, \ulb, \pl)&\;\;\text{for}\;\;y<0
\end{cases}
$,
so we have
\begin{equation*}
\lambda_0|\bm{\xi}|^2\le \sum_{i,j=1}^2a_{ij}^{(\vphil)}(x,y)\xi_i\xi_j\le \frac{1}{\lambda_0}|\bm \xi|^2
\end{equation*}
for all $(x,y)\in \N$ and $\bm{\xi}=(\xi_1, \xi_2)\in \R^2$ where $\lambda_0>0$ is a constant depending on $\Ul$ and $\gam$. Since $\bm A(x,y,\bm q)$ is smooth with respect to $\bm q$ near $D\vphil$, one can choose $\eps_1$ sufficiently small depending on $\Ul$ and $\gam$  so that if $M\sigma\le \eps_1$ in the definition of $\mcl{K}_M$, then we obtain (i) of Lemma \ref{lemma-2} for $\lambda=\frac{\lambda_0}{10}$. (iii) can be easily checked from \eqref{2-a} and \eqref{2-a3}.
\end{proof}
\end{lemma}

\begin{remark}
Lemma \ref{lemma-2}(ii) is a necessary condition to apply the result of \cite{Li}. See Theorem 1.1 of \cite{Li} for details.
\end{remark}

\begin{proposition}
\label{lemma-4}
There exists a constant $C$ depending only on $\Ul, \gam$ and $\alp$ so that if $M\sigma\le \eps_1$ for $\eps_1$ in Lemma \ref{lemma-2} and $R\ge 10$, then for any $\phi\in \mcl{K}_M$, the boundary value problem of \eqref{2-a1} and \eqref{2-a2} has unique weak solution $\psi\in H^1(\NR)\cap C^0(\ol{\NR})$ satisfying the estimate
\begin{equation}
\label{3-d5}
\begin{split}
&\|\psi\|_{2,\alp,\NhRp}^{(-1-\alp,\Gamp\cup\Gamz)}+\|\psi\|_{2,\alp,\NhRm}^{(-1-\alp,\Gamm\cup\Gamz)}\\
&\phantom{aaa}\le C(\|\Fz\|_{1,\alp,\N^+}+\|\Fz\|_{1,\alp,\N^-}+\|g\|_{1,\alp,\N}+\|Dg\|_{L^2(\N)})
\end{split}
\end{equation}
where we set
$
\N_R^+:=\{(x,y)\in \N_R: y>0\},\quad \N_R^-:=\{(x,y)\in \N_R: y<0\}.
$
\begin{proof}
\emph{(Step 1)}
Since the equation \eqref{2-a1} is uniformly elliptic, the boundary value problem of \eqref{2-a1} and \eqref{2-a2} has unique weak solution $\psi$ in $H^1(\NR)$. Also, the Harnack inequality(see \cite[Theorem 4.17, Corollary 4.18]{Ha-L}) implies that $\psi$ is continuous across $\Gamz\cap\NR$.

By Proposition 3.2 of \cite{Li}, we can choose $r_1(<R/10)$ small and $C$ depending only on $\Ul, \gam$ and $\alp$ such that for any $Z_0\in \N_{R}\cap \{(x,y):|x|<R-2r_1, |y|<\frac{r_1}{10}\}$, there is a continuous piecewise linear function $L^{(Z_0)}$ satisfying
\begin{align}
\label{3-d4}
&|\psi({\bf x})-L^{(Z_0)}({\bf x})|\le C\mu_0|{\bf x}-Z_0|^{1+\alp}\quad\text{for all}\;\; {\bf x}\in B_{r_1}(Z_0),\\
\label{3-f1}
&L^{(Z_0)}\in C^1(\ol{\N^+})\cap C^1(\ol{\N^-}),\quad
\|DL^{(Z_0)}\|_{L^{\infty}(\R^2)}\le C\mu_0
\end{align}
with
\begin{equation}
\label{mu}
\mu_0=\|\psi\|_{L^{\infty}(\N_R)}+\|\Fz\|_{1,\alp,\N^+}+\|\Fz\|_{1,\alp,\N^-}+\|g\|_{1,\alp,\N}.
\end{equation}
Then, by using the method of the proof for Theorem 1.1 in \cite{Li}, we obtain
\begin{equation}
\label{est1}
\|\psi\|_{1,\alp,\N^+_{4R/5}\cap \{|y|<\frac{r_1}{10}\}}+\|\psi\|_{1,\alp,\N^-_{4R/5}\cap \{|y|<\frac{r_1}{10}\}}\le C\mu_0
\end{equation}
for a constant $C$ depending only on $\alp$. Differently from Theorem 1.1 of \cite{Li}, we note that the regularity of $\psi$ is not weaker than $C^{1,\alp}$ because $\Gamz$ is flat thus in $C^{\infty}$. By Theorem 8.33 in \cite{GilbargTrudinger} and Lemma \ref{lemma-2}, we also have $\|\psi\|_{1,\alp,\N^+_{4R/5}\cap\{|y|>\frac{r_1}{20}\}}+\|\psi\|_{1,\alp,\N^-_{4R/5}\cap\{|y|>\frac{r_1}{20}\}}\le C\mu_0 $ for a constant $C$ depending only on $\Ul, \gam$ and $\alp$, then combining this estimate with \eqref{est1} yields
\begin{equation}
\label{est2}
\|\psi\|_{1,\alp,\N^+_{4R/5}}+\|\psi\|_{1,\alp,\N^-_{4R/5}}\le C\mu_0.
\end{equation}

For a fixed point $X_0\in \N^+_{R/2}$, let $2d_*:=dist(X_0, \Gamp\cup\Gamz)$. It suffices to consider the case of $2d_*<r_1$. Set
$
W^{(X_0)}(\eta):=\frac{\psi(X_0+d_*\eta)-\psi(X_0)}{d_*^{1+\alp}}
$
for $\eta\in B_1^{(X_0)}:=\{\eta\in B_1(0): X_0+d_*\eta\in \ol{\N^+_{R-4r_1}}\}$. Then the standard elliptic interior estimates yield
\begin{equation}
\label{3-g3}
\|W^{(X_0)}\|_{2,\alp, B_{1/2}^{(X_0)}}\le C(\|\psi\|_{1,\alp,\N^+_{3R/4}}+\|\Fz\|_{1,\alp,\N^+}),
\end{equation}
and by scaling back and combining \eqref{3-g3} with \eqref{est2}, we get
$
\|\psi\|_{2,\alp,\N^+_{R/2}}^{(-1-\alp, \Gamp\cup\Gamz)}\le C\mu_0.$
Repeating the same argument for points in $\N^-_{R/2}$, we get
\begin{equation}
\label{est3}
\|\psi\|_{2,\alp,\NhRp}^{(-1-\alp,\Gamp\cup\Gamz)}+\|\psi\|_{2,\alp,\NhRm}^{(-1-\alp,\Gamm\cup\Gamz)}\le C\mu_0.
\end{equation}
\emph{(Step 2)}
In order to finish the proof, it remains to estimate $\|\psi\|_{L^{\infty}(\N_R)}$. Set $u:=\psi-g$ then $u\in H^1_0(\N_R)$ satisfies
\begin{equation}
\label{1-b2}
\int_{\N_R}a^{(\phi)}_{ij}\der_ju\der_i\zeta=\int_{\N_R}\Fz\cdot D\zeta-a^{(\phi)}_{ij}\der_j g\der_i\zeta\;\;\tx{for any}\;\;\zeta\in H^1_0(\N_R).
\end{equation}
By \eqref{F1}, we have
\begin{equation*}
\int_{\N_R}\Fz\cdot D\zeta
=\int_{\N_R\cap\{|x|<2\}}\Fz\cdot D\zeta+(\pl-\pr)\int_{\N_R\cap\{x\ge 2\}}\der_y \zeta=\int_{\N_R\cap\{|x|<2\}}\Fz\cdot D\zeta.
\end{equation*}
Plug $\zeta=u$ into \eqref{1-b2}. Then, by Lemma \ref{lemma-2}(i) and H\"{o}lder inequality, we get
\begin{equation}
\label{3-b1}
\int_{\N_R}|Du|^2\le C(\|\Fz\|^2_{L^{\infty}(\N)}+\|Dg\|^2_{L^2(\N)})=:\mu_1
\end{equation}
where $C$ depends on $\Ul$ and $\gam$ by Lemma \ref{lemma-2} but independent of $R$. Fix a positive constant $l_0$ for $\frac{1}{10}\le l_0\le 1$. For each $z_0\in \Gamz\cap\N_R$, set $Q_{l_0}(z_0):=\N_R\cap\{|z-x_0|<l_0\}$. Since $u=0$ on $\der\N_R$, the Poincar\'{e} inequality provides $\|u\|_{L^2(Q_{l_0}(z_0))}\le C\mu_1$ for all $z_0\in \N_R$ where $C$ depends on $l_0$ but independent of $R$. Then, by the method of Moser iteration, we can find a constant $\til{C}$ depending only on $\Ul$ and $\gam$ so that $u$ satisfies
\begin{equation}
\label{est4}
\|u\|_{L^{\infty}(\N_R)}\le \til C\mu_1.
\end{equation}
Combining \eqref{est4} with \eqref{est3}, we finally obtain \eqref{3-d5}.
\end{proof}
\end{proposition}

\subsection{Proof of Proposition \ref{proposition-1}}
Now we can prove Proposition \ref{proposition-1} easily.

We choose $\eps_0=\eps_1$ for $\eps_1$ in Lemma \ref{lemma-2}.
For each $m\in \mathbb{N}$, let $\psi^{(m)}$ be unique weak solution of \eqref{2-a1} and \eqref{2-a2} in $\N_{m+20}$ with satisfying \eqref{3-d5}. Then we can extract a subsequence, still written as $\{\psi^{(m)}\}_{m\in \mathbb{N}}$, so that the subsequence converges to a function $\psi^*\in H^1_{loc}(\N)$ in the following sense: for any $R>0$
 \begin{itemize}
 \item[(i)] $\psi^{(m)}$ uniformly converges to $\psi^*$ in $\ol\N_R$;
 \item[(ii)] $\psi^{(m)}$ converges to $\psi^*$ in $C^1$ in $\ol{\N^+_R}$ and $\ol{\N^-_R}$;
 \item[(iii)] $\psi^{(m)}$ converges to $\psi^*$ in $C^2$ in $K^+$ and $K^-$ for any  $K^+\subset\subset \N^+$, $K^-\subset\subset \N^-$.
 \end{itemize}
Also, $\psi^*$ satisfies the estimate
\begin{equation}
\label{3-d9}
\begin{split}
&\|\psi^*\|_{2,\alp,\N^+}^{(-1-\alp,\Gamp\cup\Gamz)}+\|\psi^*\|_{2,\alp,\N^-}^{(-1-\alp,\Gamm\cup\Gamz)}\\
&\phantom{aaa}\le C(\|\Fz\|_{1,\alp,\N^+}+\|\Fz\|_{1,\alp,\N^-}+\|g\|_{1,\alp,\N}+\|Dg\|_{L^2(\N)})
\end{split}
\end{equation}
for $C$ same as in \eqref{3-d5}. We claim that $\psi^*$ is the unique weak solution of \eqref{3-a4}, \eqref{3-a5} satisfying the estimates \eqref{3-a6}, \eqref{3-a7} in Proposition \ref{proposition-1}.

For each constant $R>10$, let $\chi_R$ be a smooth function satisfying
\begin{equation}
\label{3-e5}
\begin{split}
&\chi_R(x,y)=\begin{cases}
1&\text{if}\;\;|x|\le R-1\\
0&\text{if}\;\;|x|\ge R-\frac 12,
\end{cases},
\quad |D\chi_R|\le 10.
\end{split}
\end{equation}
Since each $\psi^{(m)}$ is a weak solution of \eqref{2-a1} and \eqref{2-a2}, by the dominated convergence theorem and \eqref{1-b2}, $u^*:=\psi^*-g\in H^1_{loc}(\N_R)$ satisfies
\begin{equation}
\label{3-d7}
\int_{\N} a^{(\phi)}_{ij}\der_j u^*\der_i( u^*\chi_R^2)
=\int_{\N}\Fz\cdot D(u^*\chi_R^2)-a^{(\phi)}_{ij}\der_j g\der_i(u^*\chi_R^2)
\end{equation}
from which we get
\begin{equation}
\label{1-g3}
\int_{\N\cap \{|x|<R-1\}}|D\psi^*|^2\;dX
\le C\mu_1
\end{equation}
for $\mu_1$ in \eqref{3-b1} by \eqref{est4}
where the constant $C$ depends on $\Ul$ and $\gam$ but independent of $R$. Since $R$ can be arbitrarily large, we get
\begin{equation}
\label{3-e1}
\lim_{|Q|\to \infty}\int_{\N\cap\{|x-Q|\le 2\}}|D\psi^*|^2\;dx=0.
\end{equation}
Since $g_{\pm}\in H^1(\R)$,  we have
$
\underset{Q\to\infty}{\lim}\int_{|x-Q|\le 2}|g_+(x,m_0^+)|^2dx=\underset{Q\to\infty}{\lim}\int_{|x-Q|\le 2}|g_-(x,m_0^-)|^2dx=0.
$
Then, expressing $\psi^*(x,y)$ as
$
\psi^*(x,y)=\begin{cases}
g(x,m_0^+)-\int_y^{m_0^+}\psi^*_y(x,t)dt&\text{for}\;\;y>0\\
g(x,-m_0^-)+\int_{-m_0^-}^y\psi^*_y(x,t)dt&\text{for}\;\;y<0
\end{cases},
$
we can easily show that
\begin{equation}
\label{3-e2}
\lim_{|Q|\to \infty}\int_{|x-Q|\le 2}(\psi^*)^2\;dX=0,
\end{equation}
and this implies that
\begin{equation}
\label{3-e4}
\lim_{|Q|\to \infty}\|\psi^*\|_{L^{\infty}(|x-Q|\le \frac 32)}=0.
\end{equation}
by \eqref{F1} and \cite[Theorem 8.17 and 8.25]{GilbargTrudinger}. Repeating the argument of \emph{(Step 1)} in the proof of Proposition \ref{lemma-4} with using \eqref{1-a} and \eqref{F1}, we can show that $\psi^*$ satisfies \eqref{3-a7}.

Suppose that $\psi_1$ and $\psi_2$ are weak solutions to \eqref{3-a4} and \eqref{3-a5}, and that they satisfy \eqref{3-a6} and \eqref{3-a7}. Then, $\til u:=\psi_1-\psi_2$ satisfies
$
\int_{\N} a_{ij}^{(\phi)}\der_j \til u\der_i(\til u\chi_R^2)dX=0
$ for any $R\ge 10$,
and from this and \eqref{3-a7}, we get
\begin{equation*}
\int_{\N}|D\til u|^2\le C\lim_{R\to \infty}\int_{\N\cap\{R-1<|x|<R-\frac 12\}}|\til u|^2=0.
\end{equation*}
This implies $Du= 0$ in $\N^{\pm}$, therefore we get $u=k$ in $\N$ for some constants $k$. Since $\psi_1=\psi_2$ on $\Gam^{\pm}$, we conclude that $k=0$ thus $\psi_1=\psi_2$ in $\N$. The proof of Proposition \ref{proposition-1} is complete. $\Box$

\subsection{Proof of Theorem \ref{theorem-2}}
\label{subsec-3-2}
Finally, we prove Theorem \ref{theorem-2}.

Fix $\om_{\pm}$ and $h_{\pm}$ with satisfying \eqref{1-a}--\eqref{2-a6} and \eqref{2-f2}.
By Proposition \ref{proposition-1}, if $M\sigma\le \eps_1$ for $\eps_1$ in Lemma \ref{lemma-2}, then for any $\phi\in \mcl{K}_M$, the linear boundary value problem of \eqref{3-a4} and \eqref{3-a5} associated with $\phi$ has  unique weak solution $\psi^{(\phi)}\in H^1_{loc}(\N)$ satisfying the estimates \eqref{3-a6} and \eqref{3-a7}. We define a mapping $\mcl{I}$ by
\begin{equation*}
\mcl{I}: \phi\mapsto \psi^{(\phi)}+\vphiz
\end{equation*}
for $\vphi_0$ given by \eqref{3-b7}. By \eqref{1-a}, \eqref{2-a4}, \eqref{2-f2}, \eqref{2-f3}, \eqref{3-b7}, \eqref{F}, \eqref{3-a5} and \eqref{3-a6}, we have
\begin{equation*}
\|\mcl{I}(\phi)-\vphi_0\|_{2,\alp,\N^+}^{(-1-\alp,\Gamp\cup\Gamz)}+\|\mcl{I}(\phi)-\vphi_0\|_{2,\alp,\N^-}^{(-1-\alp,\Gamm\cup\Gamz)}\le C^{\flat}\sigma
\end{equation*}
for a constant $C^{\flat}$ depending only on $\Ul, \gam$ and $\alp$.
We choose $M$ and $\sigma_1$ by
\begin{equation}
\label{3-f9}
M=\max\{8C^{\flat},2\},\quad \sigma_1=\min\{\frac{\eps_1}{2M}, \frac{\sigma^{\sharp}}{2}, 1\}
\end{equation}
for $\sigma^{\sharp}$ in \eqref{1-g7} and $\eps_1$ in Lemma \ref{lemma-2}.
For such choices of $M$ and $\sigma_1$, the iteration mapping $\mcl{I}$ maps $\mcl{K}_M$ into itself wherever $\sigma\le \sigma_1$ for $\sigma$ from \eqref{1-a} and \eqref{2-a4}. We point out that the choices of $M$ and $\sigma_1$ in \eqref{3-f9} depend only on $\Ul, \gam$ and $\alp$.

We claim that $\mcl{I}:\mcl{K}_M\to \mcl{K}_M$ is continuous in $C^{1,\alp/2}_{(-1)}(\ol{\N^+})\cap C^{1,\alp/2}_{(-1)}(\ol{\N^-})$(see Section 5.1 in \cite{Ch-F3} for the definition of $C^{1,\alp/2}_{(-1)}$).
Suppose that a sequence $\{\phi_k\}$ in $\mcl{K}_M$ converges to $\bar{\phi}\in \mcl{K}_M$ in $C^{1,\alp/2}_{(-1)}(\ol{\N^+})\cap C^{1,\alp/2}_{(-1)}(\ol{\N^-})$. Let us set $\psi_k:=\mcl{I}(\phi_k)-\vphi_0$ for each $k\in \mathbb{N}$ and $\psi:=\mcl{I}(\phi)-\vphi_0$.
By \eqref{3-a6}, any subsequence of $\{\psi_k\}$ has its own subsequence that converges to a function $\psi^*$ with $\psi^*+\vphi_0\in \mcl{K}_M$ in $C^{1,\alp/2}_{(-1)}(\ol{\N^+})\cap C^{1,\alp/2}_{(-1)}(\ol{\N^-})$, and such $\psi^*$  is a weak solution to
\begin{equation*}
\sum_{i,j=1}^2\der_i(a^{(\phi)}_{ij}\der_j\psi^*)=div \Fz\quad\text{in}\;\;\N,\quad \psi^*=g\quad\text{on}\;\;\Gamp\cup\Gamm.
\end{equation*}
By repeating the argument in the proof of Proposition \ref{proposition-1}, we can show that $\psi^*=\bar{\psi}$ in $\N$. This implies that $\psi_k$ converges to $\bar{\psi}$ in $C^{1,\alp/2}_{(-1)}(\ol{N^+})\cap C^{1,\alp/2}_{(-1)}(\ol{N^-})$. Hence, $\mcl{I}:\mcl{K}_M\to\mcl{K}_M$ is continuous in $C^{1,\alp/2}_{(-1)}(\ol{N^+})\cap C^{1,\alp/2}_{(-1)}(\ol{N^-})$.

As pointed out earlier, $\mcl{K}_M$ is a convex and compact subset of $C^{1,\alp/2}_{(-1)}(\ol{\N^+})\cap C^{1,\alp/2}_{(-1)}(\ol{\N^-})$. Then, by the Schauder fixed point theorem, we conclude that for any given $h_{\pm}$ satisfying \eqref{1-a}--\eqref{2-a6}, $\mcl{I}$ has a fixed point $\vphi^{\sharp}$ in $\mcl{K}_M$. By Lemma \ref{lemma-2}, $\vphi^{\sharp}$ satisfies (ii) of Theorem \ref{theorem-2}. Also by Proposition \ref{proposition-1} and \eqref{3-b7}, $\vphi^{\sharp}$ is a weak solution of \eqref{4-6} in $\N$, and satisfies the equation \eqref{4-6} pointwisely in $\N^{\pm}$. Then, for any $\xi\in C^{\infty}_{0}(\N)$, we have
\begin{equation*}
\int_{\Gamz}[A_2(x,y,D\vphi^{\sharp})]_{\Gamz}\xi=\int \bm A(x,y,D\vphi^{\sharp})\cdot D\xi
+\int_{\N^+\cup\N^-} \xi div \bm A(x,y,D\vphi^{\sharp})=0,
\end{equation*}
so $\vphi^{\sharp}$ satisfies the R-H condition \eqref{1-16}. We may reduce $\sigma_1$ in \eqref{3-f9} further so that \eqref{2-e4} implies $(v)$ of Theorem \ref{theorem-2}. Then $\vphi^{\sharp}$ satisfies all the properties stated in Theorem \ref{theorem-2}.

Given $h=(h_+,h_-)$, let $\vphi^{(1)}$ and $\vphi^{(2)}$ be weak solutions of the boundary value problem of \eqref{4-6} and \eqref{k} with $k_{\pm}=0$ where $\vphi^{(1)}$ and $\vphi^{(2)}$ satisfy all the properties stated in Theorem \ref{theorem-2} as well. Let us set $\psi^{(j)}:=\vphi^{(j)}-\vphi_0$ for $j=1,2$. Then $\psi^{(1)}-\psi^{(2)}$ satisfies
\begin{equation}
\label{3-e6}
\begin{split}
&\sum_{i,j=1}^2\der_i
\bigl(a_{ij}^{(\vphi^{(1)})}\der_j(\psi^{(1)}-\psi^{(2)})
+(a_{ij}^{(\vphi_{(1)})}-a_{ij}^{(\vphi^{(2)})})\der_j\psi^{(2)}\bigr)=0
\quad\text{in}\;\;\N,\\
&\psi^{(1)}-\psi^{(2)}=0\;\;\tx{on}\;\;\der\N.
\end{split}
\end{equation}
By \eqref{4-7} and \eqref{2-a}, we can rewrite \eqref{3-e6} as
\begin{equation*}
\sum_{i,j=1}^2\der_i\bigl((a_{ij}^{(\vphi^{(1)})}+b_{ij})\der_j(\psi^{(1)}-\psi^{(2)})\bigr)=0
\end{equation*}
with $\|b_{ij}\|_{C^{\alp}(\ol{\N^{\pm}})}\le CM\sigma$. We again reduce $\sigma_1>0$ depending on $\Ul$ and $\gam$ to have
\begin{equation*}
\frac{\lambda'}{2}|\bm{\eta}|^2\le \sum_{i,j=1}^2(a_{ij}^{(\vphi^{(1)})}+b_{ij})\eta_i\eta_j\le \frac{2}{\lambda'}|\bm{\eta}|^2\quad\text{in}\;\;\N\;\;\tx{for all}\;\;\bm{\eta}=(\eta_1, \eta_2)\in \R^2
\end{equation*}
for some constant $\lambda'>0$. Then, repeating the argument in the proof of Proposition \ref{proposition-1}, we conclude that
$
\psi^{(1)}=\psi^{(2)}
$
in $\N$,
thus
$
\vphi^{(1)}=\phi^{(2)}
$ in $\N$. Finally, we choose $C_2=2M$, then the proof of Theorem \ref{theorem-2} is complete. $\Box$

\section{Proof of Theorem \ref{theorem-1} and Theorem \ref{theorem-1-2}}

\subsection{Proof of Theorem \ref{theorem-1}}
First, we choose $\sigma_0=\sigma_1$ for $\sigma_1$ in Theorem \ref{theorem-2}, and
fix $\om_{\pm}$ satisfying \eqref{2-f2} and $h_{\pm}$ satisfying \eqref{1-a}--\eqref{2-a6} for $\sigma\in (0,\sigma_0]$. Let $\vphi$ be the corresponding solution satisfying all the properties stated in Theorem \ref{theorem-2}. Let $(u, v, p)$ be given by \eqref{1-b3} from $\vphi$ with $\rho$ determined by \eqref{2-e3}, then $V=(u,v,p,\rho)$ satisfies all the properties stated in Problem \ref{problem3} as well as the estimates
\begin{equation}
\label{4-a1}
\begin{split}
&\|V-\Ul^+\|_{1,\alp,\N^+}^{(-\alp,\Gamp\cup\Gamz)}+\|V-\Ul^-\|_{1,\alp,\N^-}^{(-\alp,\Gamm\cup\Gamz)}
\le C\sigma,\\
&\lim_{R\to 0}\|V-\Ul^+\|_{1,\alp,\N^+\setminus\{x\ge -R\}}^{(-\alp,\Gamp\cup\Gamz)}
=\lim_{R\to 0}\|V-\Ul^-\|_{1,\alp,\N^-\setminus\{x\ge -R\}}^{(-\alp,\Gamm\cup\Gamz)}=0,\\
&\lim_{R\to 0}\|V-\Ur^+\|_{1,\alp,\N^+\setminus\{x\le R\}}^{(-\alp,\Gamp\cup\Gamz)}
=\lim_{R\to 0}\|V-\Ur^-\|_{1,\alp,\N^-\setminus\{x\le R\}}^{(-\alp,\Gamm\cup\Gamz)}=0.
\end{split}
\end{equation}
By \eqref{4-a1}, the inverse Euler-Lagrange transformation $\mathfrak{T}$ given by \eqref{1-15} is continuously differentiable in $\N^{\pm}$ and invertible with the estimates
$\|\mathfrak{T}\|_{C^{1,\alp}(\ol{\N^{\pm}})}\le C$ and $\|\mathfrak{T}^{-1}\|_{C^{1,\alp}(\ol{\til{\N}^{\pm}})}\le C$. From this, one can directly check that
$
\frac 1C d_L^+(X)\le d^+_E(X)\le Cd_L^+(X)\;\;\tx{for any}\;\;X\in \N^+
$
where we set
$
d^+_L(X):=dist(X, \Gamp\cup\Gamz)$ and $d^+_E(X):=dist(\mathfrak T(X), \mathfrak T(\Gamp\cup \Gamz))$, that is, $d^+_L$ is a distance function in the Lagrangian coordinates, and $d^+_E$ is a distance function in the Eulerian coordinates.
Similarly, if we set $d^-_L(X):=dist(X, \Gamm\cup\Gamz)$ and $d^-_E(X):=dist(\mathfrak{T}(X), \mathfrak{T}(\Gamm\cup \Gamz))$ for $X\in \N^-$, then we also have
$
\frac 1C d^-_L(X)\le d^-_E(X)\le Cd^-_L(X).
$

If we set
$$
\cd:=\mathfrak{T}(\Gamz)\;\;\tx{and}\;\;U:=V\circ \mathfrak{T}^{-1},
$$
then $U$ satisfies all the properties of Problem \ref{problem1} as well as (iii) and (iv) of Theorem \ref{theorem-1} where we choose $C_0$ as $CC_1$ for some constant $C$ depending only on $\Ul, \gam$ and $\alp$.
By the definition of $\mathfrak{T}$ in \eqref{1-15}, the contact discontinuity $\cd$ is given by
\begin{equation}
\label{1-g1}
\begin{split}
\cd&=\{(x,\gd(x)):\gd(x)=\int_{-m_0^-}^0\frac{1}{\rho u(x,t)}dt+h_-(x), x\in \R\}\\
&\tx{or}=\{(x,\gd(x)):\gd(x)=h_+(x)-\int_{0}^{m_0^+}\frac{1}{\rho u(x,t)}dt, x\in \R\}.
\end{split}
\end{equation}
By \eqref{4-5}, \eqref{k} with $k_{\pm}=0$ and \eqref{2-f1}, we can express $\gd$ as
\begin{equation}
\label{1-g6}
\gd(x)
=\int_{-m_0^-}^0\vphi_y(x,t)dt+h_-(x)=\vphi(x,0)=(\vphi-\vphil)(x,0).
\end{equation}
Then, (i) and (ii) of Theorem \ref{theorem-1} follow from \eqref{2-f1}, \eqref{2-e4} and \eqref{2-e5}. Theorem \ref{theorem-1}(v) follows from Theorem \ref{theorem-2}(v) and \eqref{4-5}.

Now it remains to verify the uniqueness in Theorem \ref{theorem-1}.
For fixed $\om_{\pm}$ and $h_{\pm}$, let $U^{(1)}$ and $U^{(2)}$ be two solutions of Problem \ref{problem1} with satisfying all the properties stated in Theorem \ref{theorem-1}, and let $T^{(j)}$ for $j=1,2$ be defined from $U^{(j)}$ by \eqref{1-1}, and set $V^{(j)}:=U^{(j)}\circ (T^{(j)})^{-1}$ in $\N$. Then, each $V^{(j)}=(u^{(j)}, v^{(j)}, p^{(j)}, \rho^{(j)})$ satisfies \eqref{4-a1}. For $j=1,2$, if we can find $\vphi^{(j)}$ satisfying
\begin{align}
\label{4-a5}
&(\vphi_x^{(j)},\vphi^{(j)}_y)=(\frac{v^{(j)}}{u^{(j)}},\frac{1}{\rho^{(j)}u^{(j)}})\quad\text{in}\;\;\N^{\pm},\\
\label{4-a6}
&\vphi^{(j)}=h_{\pm}\;\;\text{on}\;\;\Gam^{\pm}
\end{align}
then $\vphi^{(j)}$ satisfies all the properties in Theorem \ref{theorem-2} by reducing $\sigma_0>0$ if necessary. Then, following the proof of Theorem \ref{theorem-2}, we can easily show that $\vphi^{(1)}=\vphi^{(2)}$ in $\N$ so $U^{(1)}=U^{(2)}$ in $\til{\N}$. The proof is complete.
\hfill $\Box$

\subsection{Proof of Theorem \ref{theorem-1-2}} For $\om_{\pm}$ and $h_{\pm}$ fixed, let $(U,\gd)$ be the corresponding solution satisfying all the properties in Theorem \ref{theorem-1}, and let $\vphi^{\natural}$ be the solution corresponding to $(U,\gd)$ so that $\vphi^{\natural}$ satisfies all the properties in Theorem \ref{theorem-2}. Then, by \eqref{4-5} and \eqref{1-g1}, the contact discontinuity function $\gd$ is given by
\begin{equation*}
\gd(x)=h_+(x)-\int_0^{m_0^+}\vphi_y(x,t) dt.
\end{equation*}
Then, by \eqref{2-a8}, \eqref{4-5} and \eqref{3-a5}, we have
$
\gd(x)-\om_*\eta(x)=g_+(x)-\int_0^{m_0^+}\der_y(\vphi-\vphiz)(x,t)dt
$
for $\omega_*$ in Lemma \ref{lemma-2-1} and $\eta$ defined by \eqref{1-g2},
and this implies
\begin{equation*}
\|\gd-\om_*\eta\|_{L^2(\R)}\le 2(\|g_+\|_{L^2(\R)}+m_0^+\|D(\vphi-\vphiz)\|_{L^2(\N^+)}).
\end{equation*}
Then, by \eqref{1-a}, \eqref{2-e7}, \eqref{F1} and \eqref{1-g3}, we obtain \eqref{1-g5}.

\eqref{1-g4} can be similarly proved by using \eqref{2-f3}, \eqref{2-b9}, \eqref{1-g3} and Theorem \ref{theorem-1}(v).
\hfill $\Box$

\bigskip
{\bf Acknowledgments.} The author thanks Gui-Qiang Chen at University of Oxford for motivating to work on this problem, and Mikhail Feldman at University of Wisconsin-Madison for helpful discussion.

\end{document}